\documentclass[10pt]{article}
\title{Representation theory of Yang-Mills algebras}
\author{Estanislao Herscovich and Andrea Solotar
\thanks{This work has been supported by the projects UBACYTX169 and X212, CONICET-CNRS, PICS 3410 and PIP-CONICET 5099.
The first author is a CONICET fellow.
The second author is a research member of CONICET (Argentina) and a
Regular Associate of ICTP Associate Scheme.}}
\date{}

\input{xy}
\xyoption{all}
\xyoption{poly}
\usepackage[all]{xy}
\usepackage{amsmath,amsthm,amsfonts,amssymb,stmaryrd}
\usepackage{latexsym}
\usepackage{color,palatino}
\usepackage{graphicx}
\usepackage{float}  
\usepackage{fullpage}
\usepackage[small, bf, margin=90pt, tableposition=bottom]{caption}
\RequirePackage{amssymb}
\RequirePackage[T1]{fontenc}
\usepackage[letterpaper]{geometry}

\usepackage[
  breaklinks,
  naturalnames,
  ]{hyperref}

\newtheorem{teo}{Theorem}[section]
\newtheorem{coro}[teo]{Corollary}
\newtheorem{lema}[teo]{Lemma}
\newtheorem{prop}[teo]{Proposition}
\newtheorem{defi}[teo]{Definition}
\newtheorem{obs}[teo]{Remark}
\newtheorem{rem}[teo]{Remark}

\newtheorem{ejem}[teo]{Example}

\numberwithin{equation}{section}                    
\def\id{{\mathrm{id}}}
\let\oldqed\qed
\renewcommand\qed{\oldqed\par\bigskip}

\newcommand\cl[1]{{\langle#1\rangle}}

\newcommand\ZZ{{\mathbb{Z}}}
\newcommand\NN{{\mathbb{N}}}

\newcommand\gl{{\mathfrak{gl}}}

\def\Aut{{\mathrm {Aut}}}
\def\Prim{{\mathrm{Prim}}}

\def\A{{\mathcal A}}
\def\B{{\mathcal B}}
\def\C{{\mathcal C}}
\def\D{{\mathcal D}}
\def\E{{\mathcal E}}
\def\F{{\mathcal F}}
\def\G{{\mathcal G}}
\def\H{{\mathcal H}}
\def\I{{\mathcal I}}

\def\L{{\mathcal L}}

\def\O{{\mathcal O}}
\def\R{{\mathcal R}}

\def\U{{\mathcal U}}
\def\W{{\mathcal W}}
\def\Z{{\mathcal Z}}

\def\YM{{\mathrm {YM}}}
\def\TYM{{\mathrm {TYM}}}

\def\ad{{\mathrm {ad}}}

\def\Ker{{\mathrm {Ker}}}

\def\Mod{{\mathrm {Mod}}}
\def\mod{{\mathrm {mod}}}
\def\Frac{{\mathrm {Frac}}}
\def\Ind{{\mathrm {Ind}}}
\def\End{{\mathrm {End}}}
\def\tr{{\mathrm {tr}}}
\def\Alg{{\mathrm {Alg}}}

\def\LieAlg{{\mathrm {LieAlg}}}

\def\Der{{\mathrm {Der}}}
\def\InnDer{{\mathrm {InnDer}}}

\def\Tor{{\mathrm {Tor}}}
\def\Ind{{\mathrm {Ind}}}
\def\GW{{\mathcal {GW}}}

\def\Alg{{\mathrm {Alg}}}

\def\g{{\mathfrak g}}
\def\h{{\mathfrak h}}
\def\n{{\mathfrak n}}
\def\ym{{\mathfrak{ym}}}
\def\f{{\mathfrak f}}

\def\tym{{\mathfrak{tym}}}

\def\place{{-}}

\def\MyNode{\ifcase\xypolynode\or
      (W \otimes X) \otimes (Y \otimes Z)
    \or
      ((W \otimes X) \otimes Y) \otimes Z
    \or
      (W \otimes (X \otimes Y)) \otimes Z
    \or
      W \otimes ((X \otimes Y) \otimes Z)
    \or
      W \otimes (X \otimes (Y \otimes Z))
    \fi
  }%

\def\MyNodes{\ifcase\xypolynode
    \or
      X \otimes Y
    \or
      X \otimes (e \otimes Y)
    \or
     (X \otimes e) \otimes Y
    \fi
  }%

\def\MyNodessr{\ifcase\xypolynode
    \or
      X \otimes Y
    \or
      e \otimes (X \otimes Y)
    \or
     (e \otimes X) \otimes Y
    \fi
  }%

\def\MyNodessl{\ifcase\xypolynode
    \or
      X \otimes Y
    \or
      (X \otimes Y) \otimes e
    \or
       X \otimes (Y \otimes e)
    \fi
  }%
\begin{document}

\maketitle
\begin{abstract}
   {
   The aim of this article is to describe families of representations of the Yang-Mills algebras $\YM(n)$
   ($n \in \NN_{\geq 2}$) defined by A. Connes and M. Dubois-Violette in \cite{CD1}.
   We first describe irreducible finite dimensional representations.
   Next, we provide families of infinite dimensional representations of $\YM(n)$, big enough to separate points of the algebra.
   In order to prove this result, we prove and use that all Weyl algebras $A_{r}(k)$ are epimorphic images of $\YM(n)$.
   }
\end{abstract}



\section{Introduction}

This article is devoted to the study of the representation theory of Yang-Mills algebras.
Very little is known on this subject.
Our goal is to describe families of representations which, although they do not cover the whole category of representations,
are large enough to distinguish elements of the Yang-Mills algebra.

In order to describe our results in more detail, let us recall the definition of Yang-Mills algebras by A. Connes and
M. Dubois-Violette in \cite{CD1}.
Given a positive integer $n \geq 2$, the Lie Yang-Mills algebra over an algebraically closed field $k$ of characteristic zero is
\[     \ym(n) = \f(n)/\cl{\{ \sum_{i=1}^{n} [x_{i},[x_{i},x_{j}]] : j= 1, \dots, n \}},     \]
where $\f(n)$ is the free Lie algebra in $n$ generators $x_{1}, \dots, x_{n}$.
Its associative enveloping algebra $\U(\ym(n))$ will be denoted $\YM(n)$.
It is, for each $n$, a cubic Koszul algebra of global dimension $3$ and we shall see that it is noetherian if and only if $n=2$.

They have been recently studied mainly due to its physical applications, since they arise as limits of algebras
appearing in the gauge theory of $D$-branes and open string theory in background independent formulations (cf. \cite{Ne1, Mov1}).
As pointed out by M. Douglas, the case $n=4$ appears, by means of the AdS/CFT correspondence,
as a reformulation of a string theory in the anti de Sitter space (cf. \cite{Doug1}).

In \cite{HKL1}, the authors discuss a superized version of Yang-Mills algebras. 

Our main result may be formulated as follows:
\begin{teo}
Given $n \geq 3$ and $r \geq 1$, the Weyl algebra $A_{r}(k)$ is an epimorphic image of $\YM(n)$.
\end{teo}

The key ingredient of the proof is the existence of a Lie ideal in $\ym(n)$ which is free as Lie algebra. 
This ideal has already been considered by M. Movshev in \cite{Mov1}.
It allows us to define morphisms from the Yang-Mills algebras onto the Weyl algebras,
making use for this of the Kirillov orbit method.

Once this is achieved, the categories of representations of \emph{all} Weyl algebras, that have been extensively studied by
V. Bavula and V. Bekkert in \cite{BB1}, are also representations of the Yang-Mills algebras.
Thus, we provide several families of representations of the Yang-Mills algebras.
However, an easy argument using Gelfand-Kirillov dimension shows that this construction does not provide all representations. 

\bigskip

The contents of the article are as follows.
In section \ref{sec:generalities} we recall the definition and elementary properties of Yang-Mills algebras and we also study the
subcategory of finite dimensional representations, describing completely those which are finite dimensional and irreducible. 

Section \ref{sec:tym} is devoted to the description of the Lie ideal $\tym(n)$.
We give complete proofs of the fact that it is a free Lie algebra in itself.
This involves the construction of a model of the graded associative algebra $\TYM(n) = \U(\tym(n))$,
which permits us to replace the bar complex of $\TYM(n)$ by the quasi-isomorphic bar complex of the model of $\TYM(n)$.

Finally, in section \ref{sec:main} we prove our main result, Corollary \ref{coro:yangmillsweyl},
and describe the families of representations appearing in this way.

Throughout this article $k$ will denote an algebraically closed field of characteristic zero.
Given an ordered basis $\{ v_{1},\dots, v_{n} \} \subset V$ of the $k$-vector space $V$,
$\{ v_{1}^{*}, \dots, v_{n}^{*} \} \subset V^{*}$ will denote its dual basis. 
Also, we shall identify $\g$ with its image inside $\U(\g)$ via the canonical morphism $\g \rightarrow \U(\g)$.

Given an associative or Lie $k$-algebra $A$ and a subgroup $G$ of $\ZZ$, we shall denote ${}_{A}^{G}\Mod$, $\Mod_{A}^{G}$,
${}_{A}^{G}\mod$ and $\mod_{A}^{G}$ the categories of $G$-graded left and right $A$-modules, and finite dimensional $G$-graded left
and right $A$-modules, respectively.
We will also denote ${}_{k}^{G}\Alg$ and ${}_{k}^{G}\LieAlg$ the categories of $G$-graded associative and Lie algebras, respectively.
We notice that if $G$ is trivial, each definition yields the non-graded case.

We would like to thank Jacques Alev, Jorge Vargas and Michel Dubois-Violette for useful comments and remarks. 
We are indebted to Mariano Su\'arez-\'Alvarez for a careful reading of the manuscript, suggestions and improvements. 

\section{Generalities and finite dimensional modules}
\label{sec:generalities}

In this first section we fix notations and recall some elementary properties of the Yang-Mills algebras.
We also study the category of finite dimensional representations, describing completely all irreducible finite dimensional representations.

Let $n$ be a positive integer such that $n \geq 2$ and let $\f(n)$ be the free Lie algebra with generators $\{ x_{1}, \dots, x_{n} \}$.
This Lie algebra is trivially provided with a locally finite dimensional $\NN$-grading.

Following \cite{CD1}, the quotient Lie algebra
\[     \ym(n) = \f(n)/\cl{\{ \sum_{i=1}^{n} [x_{i},[x_{i},x_{j}]] : 1 \leq j \leq n \}},     \]
is called the \emph{Yang-Mills algebra with $n$ generators}.
The $\NN$-grading of $\f(n)$ induces an $\NN$-grading of $\ym(n)$, which is also locally finite dimensional.
We denote $\ym(m)_{i}$ the $i$-th homogeneous component, so that 
\begin{equation}
\label{eq:gradym}
     \ym(n) = \bigoplus_{j \in \NN} \ym(n)_{j},
\end{equation}
and put $\ym(n)^{l} = \bigoplus_{j = 1}^{l} \ym(n)_{j}$. 
The Lie ideal
\begin{equation}
\label{eq:tym}
   \tym(n) = \bigoplus_{j \geq 2} \ym(n)_{j}
\end{equation}
will be of considerable importance in the sequel.

The universal enveloping algebra $\U(\ym(n))$ will be denoted by $\YM(n)$ and it is 
called the \emph{(associative) Yang-Mills algebra with $n$ generators}.
If $V(n) = \text{span}_{k} (\{ x_{1} , \dots , x_{n} \})$ and 
$R(n) = \mathrm{span}_{k} \cl{\{ \sum_{i=1}^{n} [x_{i},[x_{i},x_{j}]] : 1 \leq j \leq n \}} \subset V(n)^{\otimes 3}$, it turns out that
\[     \YM(n) \simeq TV(n)/\cl{R(n)}.     \]
We shall also consider the universal enveloping algebra of the Lie ideal $\tym(n)$, which will be denoted $\TYM(n)$.
Occasionally, we will omit the index $n$ in order to simplify the notation if it is clear from the context.

It is clear to see that the Yang-Mills algebra $\YM(n)$ is a domain for any $n \in \NN$, since it is the enveloping algebra
of a Lie algebra (cf. \cite{Dix1}, Corollary 2.3.9, (ii), p. 76).

The first example of Yang-Mills algebra appears when $n = 2$.
We shall see in the sequel that it is in fact quite different from the other cases.
\begin{ejem}
\label{ejem:heisenberg}
Let $n =2$.
In this case, $\ym(2)$ is isomorphic to the \emph{Heisenberg Lie algebra} $\h_{1}$, with generators $x,y,z$, and relations $[x,y] = z$,
$[x,z] = [y,z] = 0$.
The isomorphism is given by $x_{1} \mapsto x$, $x_{2} \mapsto y$.

Alternatively, $\ym(2) \simeq \n_{3}$, where $\n_{3}$ is the Lie algebra of strictly upper triangular $3 \times 3$ matrices with coefficients in $k$.
The isomorphism is now given by $x_{1} \mapsto e_{12}$, $x_{2} \mapsto e_{23}$.

We see that $\YM(2)$ is a noetherian algebra, since $\ym(2)$ is finite dimensional.
Furthermore, since $\U(\h_{1}) \simeq A(2,-1,0)$, the Yang-Mills algebra with two generators is isomorphic to a down-up algebra,
already known to be noetherian.
\end{ejem}

We shall consider two different but related gradings on $\ym(n)$.
On the one hand, the grading given by \eqref{eq:gradym} will be called the \emph{usual grading} of the Yang-Mills algebra $\ym(n)$.
On the other hand, following \cite{Mov1}, we shall also consider the \emph{special grading} of the Yang-Mills algebra $\ym(n)$,
for which it is a graded Lie algebra concentrated in even degrees with each homogeneous space $\ym(n)_{j}$ in degree $2j$.
In this case, the Lie ideal $\tym(n)$ given in \eqref{eq:tym} is also concentrated in even degrees (strictly greater than $2$).
We will see that in fact $\tym(n)$ is isomorphic (as graded Lie algebras) to the graded free Lie
algebra on a graded vector space $W(n)$, i.e., $\tym \simeq \f_{gr}(W(n))$ (cf. \cite{Mov1}, \cite{MS1} and Section \ref{sec:tym}).

These gradings of the Lie algebra $\ym(n)$ induce respectively the \emph{usual grading} and the \emph{the special grading} on the associative algebra
$\YM(n)$.
This last one corresponds to taking the graded universal enveloping algebra of the graded Lie algebra $\ym(n)$.

In order to understand the relation among the graded and non-graded cases we present the following proposition for which we omit the proof:
\begin{prop}
\label{prop:conmutacionfuntores}
The following diagram of functors, where $\mathcal{O}$ denote the corresponding forgetful functors,
\[
\xymatrix@R-20pt 
{ 
{}^{2\ZZ}_{k}\Mod 
\ar[ddd]^(0.65){\mathcal{O}}
\ar[ddr]^{T_{gr}} 
\ar[drr]^(0.35){\mathfrak{f}_{gr}} 
& 
&
\\
& 
& 
{}^{2\ZZ}_{k}\LieAlg 
\ar[ddd]^(0.65){\mathcal{O}}
\ar[ld]^{\U_{gr}}
\\
& 
{}^{2\ZZ}_{k}\Alg 
\ar[ddd]^(0.65){\mathcal{O}} 
&
\\
{}_{k}\Mod \ar[ddr]^{T} 
\ar[drr]|(0.4825){\phantom{x}}^(0.35){\mathfrak{f}} 
& 
&
\\
& 
& 
{}_{k}\LieAlg 
\ar[ld]^{\U}
\\
& 
{}_{k}\Alg 
& 
}
\]
is commutative.
\end{prop}

\begin{obs}
\label{obs:libregraduado}
The forgetful functors ${}^{2\ZZ}_{k}\LieAlg \rightarrow {}_{k}\LieAlg$ and ${}^{2\ZZ}_{k}\Alg \rightarrow {}_{k}\Alg$ preserve and reflect free objects.
\end{obs}


The Yang-Mills algebra is not nilpotent in general, since its lower central series
\[     \ym(n) = \C^{0}(\ym(n)) \supset \C^{1}(\ym(n)) \supset \dots \supset \C^{k}(\ym(n)) \supset \dots     \]
is not finite as we will show below.
However, it is \emph{residually nilpotent}, that is, $\cap_{m \in \NN} \C^{m}(\ym(n)) = 0$.
This fact is a direct consequence of the following: $\C^{m}(\ym(n))$ is included in $\oplus_{j \geq m+1} \ym(n)_{j}$, since $\ym(n)$ is graded, 
and $\ym(n)/\C^{m}(\ym(n))$ is a finite dimensional nilpotent Lie algebra for every $m \in \NN_{0}$,
since the lower central series of $\ym(n)/\C^{m}(\ym(n))$ is
\[     \ym(n)/\C^{m}(\ym(n)) = \C^{0}(\ym(n)/\C^{m}(\ym(n))) \supset \C^{1}(\ym(n))/\C^{m}(\ym(n)) \supset
       \dots \supset \C^{m}(\ym(n))/\C^{m}(\ym(n)) = 0.     \]

Let us study in detail the ideals appearing in the lower central series.
Since the ideal of $\f(n)$
\[     I(n) = \cl{\{ \sum_{i=1}^{n} [x_{i},[x_{i},x_{j}]] : 1 \leq j \leq n \}}     \]
is homogeneous, then
\[     I(n) = \bigoplus_{j \in \NN} I(n)_{j} = \bigoplus_{j \in \NN} (I(n) \cap \f(n)_{j}).     \]
We also notice that in the free Lie algebra,
\[     \C^{k}(\f(n)) = \bigoplus_{j \geq k +1} \f(n)_{j},     \]
so that
\[     \C^{k}(\ym(n)) = \C^{k}(\f(n))/ (I(n) \cap \C^{k}(\f(n))) = \bigoplus_{j \geq k +1} \f(n)_{j}/(I(n) \cap \f(n)_{j})
        = \bigoplus_{j \geq k + 1} \f(n)_{j}/I(n)_{j} = \bigoplus_{j \geq k + 1} \ym(n)_{j}.     \]
Hence, there exists a canonical $k$-linear isomorphism $j_{l} : \ym(n)/\C^{l}(\ym(n)) \rightarrow \ym(n)^{l}$.
As a consequence, we see that if $\ym(n)$ is not finite dimensional, then $\C^{k}(\ym(n)) \neq 0$ for $k \in \NN_{0}$, and $\ym(n)$ is not nilpotent.
We shall prove below that the Yang-Mills algebra $\ym(n)$ is finite dimensional if and only if $n = 2$ (cf. Remark \ref{rem:nonfinitedimensional}).

Since we are interested in studying representations of the Yang-Mills algebra we prove the following useful lemma.
\begin{lema}
\label{lema:yangmillsnulo}
For each $l \in \NN$, the surjective Lie algebra homomorphism $\pi_{l} : \ym(n) \twoheadrightarrow \ym(n)/\C^{l}(\ym(n))$,
induces a surjective algebra homomorphism $\Pi_{l} : \U(\ym(n)) \twoheadrightarrow \U(\ym(n)/\C^{l}(\ym(n)))$.
Let $K_{l} = \Ker(\Pi_{l})$.
Then, the fact that $\cap_{l \in \NN} \C^{l}(\ym(n)) = 0$ implies that 
\[     K = \bigcap_{l \in \NN} K_{l} = 0.     \]
\end{lema}
\noindent\textbf{Proof.}
The Poincar\'e-Birkhoff-Witt Theorem says that, given a Lie algebra $\g$, there exists a canonical $k$-linear isomorphism
$\gamma : S(\g) \overset{\simeq}{\rightarrow} \U(\g)$ given by symmetrization (cf. \cite{Dix1}, 2.4.5). 
We shall denote by $\epsilon_{\g} : S(\g) \rightarrow k$ the augmentation of $S(\g)$ given by the canonical projection over the field $k$.

Since the functor $S(\place)$ is a left adjoint to the forgetful functor of commutative $k$-algebras into $k$-modules, it preserves colimits.
In particular, we obtain from the decomposition $\ym(n) = \ym(n)^{l} \oplus \C^{l}(\ym(n))$ that $S(\ym(n))$ is isomorphic to
$S(\ym(n)^{l}) \otimes S(\C^{l}(\ym(n)))$, via the $k$-algebra isomorphism $t$
induced by the $k$-linear map $v + w \mapsto v \otimes 1 + 1 \otimes w$, where $v \in \ym(n)^{l}$, $w \in \C^{l}(\ym(n))$.
The inverse of $t$ is multiplication $v \otimes w \mapsto vw$.

The surjective $k$-linear map $\pi_{l} : \ym(n) \twoheadrightarrow \ym(n)/\C^{l}\ym(n))$ induces a surjective $k$-algebra homomorphism
$P_{l} : S(\ym(n)) \twoheadrightarrow S(\ym(n)/\C^{l}(\ym(n)))$.
Analogously, the $k$-linear isomorphism $j_{l} : \ym(n)/\C^{l}(\ym(n)) \rightarrow \ym(n)^{l}$ induces a $k$-algebra isomorphism
$J_{l} : S(\ym(n)/\C^{l}\ym(n))) \rightarrow S(\ym(n)^{l})$.
The composition $J_{l} \circ P_{l}$ coincides with $(\id_{\ym(n)^{l}} \otimes \epsilon_{\C^{l}(\ym(n))}) \circ t$, and hence has kernel
$t^{-1}(S(\ym(n)^{l}) \otimes S_{+}(\C^{l}(\ym(n)))) = S(\ym(n)^{l}) S_{+}(\C^{l}(\ym(n)))$.

On the other hand, the following diagram
\[
\xymatrix { S(\ym(n)) \ar[rr]^(0.4){P_{l}} \ar[d]^{\gamma} & &
S(\ym(n)/\C^{l}(\ym(n))) \ar[d]^{\gamma}
\\
\U(\ym(n)) \ar[rr]^(0.4){\Pi_{l}} & & \U(\ym(n)/\C^{l}(\ym(n))) }
\]
is commutative.

As a consequence, $K_{l} = \gamma(S(\ym(n)^{l}) S_{+}(\C^{l}(\ym(n))))$.
Whence, using the fact that $\gamma$ is bijective,
\[     K = \bigcap_{l \in \NN} K_{l} = \bigcap_{l \in \NN} \gamma(S(\ym(n)^{l}) S_{+}(\C^{l}(\ym(n))))
         = \gamma(\bigcap_{l \in \NN} S(\ym(n)^{l}) S_{+}(\C^{l}(\ym(n)))) = 0.     \]
The last equality can be proved as follows:
taking into account that
\[     \big( \bigcap_{l \in \NN} S(\ym(n)^{l}) S_{+}(\C^{l}(\ym(n))) \big) \cap S(\ym(n)^{l}) = 0,      \]
for all $l \in \NN$, and, since $S(\ym(n)) = \bigcup_{l \in \NN} S(\ym(n)^{l})$, then
\[     \bigcap_{l \in \NN} S(\ym(n)^{l}) S_{+}(\C^{l}(\ym(n))) = 0.      \]
\qed

Let
\[     \psi : \ym(n)/\C^{l}(\ym(n)) \rightarrow \gl(V),     \]
be a representation of the quotient $\ym(n)/\C^{l}(\ym(n))$.
It provides a representation of $\ym(n)$ simply by composition with the canonical projection $\pi_{l}$.
Given a morphism $f$ between two representations $V$ and $W$ of the quotient $\ym(n)/\C^{l}(\ym(n))$, it induces
a morphism between the corresponding representations of the algebra $\ym(n)$ in a functorial way.
Hence, it yields a $k$-linear functor
\[     I_{l} : {}_{\ym(n)/\C^{l}(\ym(n))}\Mod \rightarrow {}_{\ym(n)}\Mod,     \]
which also restricts to the full subcategories of finite dimensional modules, denoted by $i_{l}$.
Moreover, since the map $\pi_{l} : \ym(n) \rightarrow \ym(n)/\C^{l}(\ym(n))$ is onto, the change-of-rings functors are fully faithful.  

Analogously, given $l, m \in \NN$, such that $l \leq m$, the homomorphism of Lie algebras
$\pi_{l \leq m} : \ym(n)/\C^{m}(\ym(n)) \rightarrow \ym(n)/\C^{l}(\ym(n))$ induced by the canonical projection gives a $k$-linear functor
\[     I_{l \leq m} : {}_{\ym(n)/\C^{l}(\ym(n))}\Mod \rightarrow {}_{\ym(n)/\C^{m}(\ym(n))}\Mod,     \]
which restricts to the full subcategories of finite dimensional modules.
We shall denote this restriction by $i_{l \leq m}$.
It is clear that $l_{m \leq p} \circ I_{l \leq m} = I_{l \leq p}$ and $I_{m} \circ I_{l \leq m} = I_{l}$.

\begin{rem}
\label{rem:irred}
The functors $I_{l \leq m}$ and $I_{l}$ ($l, m \in \NN$) preserve irreducible modules.
\end{rem}

The following proposition concerning the categories of finite dimensional modules is easy to prove.
\begin{prop}
\label{prop:mod}
Let $\phi : \ym(n) \rightarrow \gl(V)$ be a nilpotent finite dimensional representation of $\ym(n)$.
Then, there exist $m \in \NN$ and a homomorphism of Lie algebras $\phi' : \ym(n)/\C^{m}(\ym(n)) \rightarrow \gl(V)$, such that
$\phi = \phi' \circ \pi_{m}$. 
The converse is also true. 
\end{prop}
\noindent\textbf{Proof.}
Since the image of $\phi$ is a nilpotent finite dimensional subalgebra of $\gl(V)$, $\Ker(\phi)$ is finite codimensional 
and there exists an $m \in \NN$ such that $\Ker(\phi) \supset \oplus_{j \geq m} \ym(n)_{j}$.
Since $\oplus_{j \geq m} \ym(n)_{j} \supset \C^{m}(\ym(n))$, we get $\Ker(\phi) \supset \C^{m}(\ym(n))$.
Therefore $\phi$ induces a morphism $\phi' : \ym(n)/\C^{m}(\ym(n)) \rightarrow \gl(V)$ satisfying $\phi = \phi' \circ \pi_{m}$. 
The converse is clear. 
\qed

\begin{rem}
The authors do not know if the previous proposition holds if one does not ask the representation to be nilpotent. 
It is not difficult to prove that every two dimensional representation of the Yang-Mills algebra is indeed nilpotent, 
but the question remains open for higher dimensions. 
\end{rem}

\begin{coro}
The category ${}_{\ym(n)}{}_{\textrm{nil}}\mod$ of nilpotent finite dimensional modules over $\ym(n)$ is the 
filtered colimit in the category of $k$-linear categories of the categories ${}_{\ym(n)/\C^{m}(\ym(n))}\mod$ of finite dimensional 
modules over $\ym(n)/\C^{m}(\ym(n))$.
\end{coro}
\noindent\textbf{Proof.}
Let $\C$ be a $k$-linear category and let $F_{l} : {}_{\ym(n)/\C^{l}(\ym(n))}\mod \rightarrow \C$ be a collection of $k$-linear functors indexed by
$l \in \NN$ satisfying that $F_{m} \circ I_{l \leq m} = F_{l}$, for $l,m \in \NN$, $l \leq m$.
We shall define a $k$-linear functor $F : {}_{\ym(n)}{}_{\textrm{nil}}\mod \rightarrow \C$.

If $M$ is a nilpotent finite dimensional representation of $\ym(n)$ given by $\phi : \ym(n) \rightarrow \gl(M)$, using Proposition \ref{prop:mod}
we see that there exists $l \in \NN$ such that $\phi$ can be factorized as $\phi_{l} \circ \pi_{l}$ with
$\phi_{l} : \ym(n)/\C^{l}(\ym(n)) \rightarrow \gl(M)$, and hence $M$ can be considered as a module over $\ym(n)/\C^{l}(\ym(n))$, denoted $M_{l}$.
We define $F(M) = F_{l} (M_{l})$.

The functor $F$ is well-defined.
Suppose that for another $m \in \NN$ there exists $\phi_{m} : \ym(n)/\C^{m}(\ym(n)) \rightarrow \gl(M)$
such that $\phi = \phi_{m} \circ \pi_{m}$, then for $m \geq l$ there is a diagram
\[
\xymatrix@R-10pt
{
&
\ym/\C^{m}(\ym)
\ar@{->>}[dddd]_{\pi_{l \leq m}}
\ar[ddrr]^{\phi_{m}}
&
&
\\
&
&
&
\\
&
&
&
\gl(M)
\\
\ym
\ar@{->>}[uuur]^{\pi_{m}}
\ar@{->>}[rd]^{\pi_{l}}
\ar[urrr]|(0.385){\phantom{x}}^{\phi}
&
&
&
\\
&
\ym/\C^{l}(\ym)
\ar[uurr]^{\phi_{l}}
&
&
}
\]

From the definitions, all faces are commutative except maybe the one implying $\phi_{l} \circ \pi_{l \leq m} = \phi_{m}$.
Since $\pi_{m}$ is surjective, the previous equality is satisfied if and only if $\phi_{l} \circ \pi_{l \leq m} \circ \pi_{m}
= \phi_{m} \circ \pi_{m}$.
But
\[     \phi_{l} \circ \pi_{l \leq m} \circ \pi_{m} = \phi_{l} \circ \pi_{l} = \phi = \phi_{m} \circ \pi_{m},     \]
so $M_{m} = I_{k \leq m} (M_{k})$ and $F_{m}(M_{m}) = F_{m} \circ I_{l \leq m} (M_{l}) = F_{l} (M_{l})$.

Let $M$ and $N$ be two nilpotent finite dimensional representations of $\ym(n)$ by means of $\phi : \ym(n) \rightarrow \gl(M)$ and
$\psi : \ym(n) \rightarrow \gl(N)$ and let $f : M \rightarrow N$ be a module homomorphism.
By Proposition \ref{prop:mod} there exists $l \in \NN$ such that both $\phi$ and $\psi$ can be factorized as $\phi_{l} \circ \pi_{l}$ and
$\psi_{l} \circ \pi_{l}$, respectively.
As before, we shall denote $M_{l}$ and $N_{l}$ these modules.
It follows directly from the definitions that $f$ is also a morphism of $\ym(n)/\C^{l}(\ym(n))$-modules, denoted $f_{l}$.
Take $F(f) = F_{l} (f_{l})$.
It is clearly well-defined.
\qed

The previous proposition says that every irreducible nilpotent finite dimensional $\ym(n)$-module is an irreducible finite dimensional
module over $\ym(n)/\C^{l}(\ym(n))$, for some $l \in \NN$. 
Since the latter is a finite dimensional nilpotent Lie algebra, it suffices to find irreducible finite dimensional modules over this
kind of Lie algebras.
It is in fact well-known that these representations are one dimensional (cf. \cite{Dix1}, Coro. 1.3.13).
The following lemma provides a description of them. 
\begin{lema}
\label{lema:nil}
Let $\g$ be a nilpotent finite dimensional Lie algebra over $k$.
Every irreducible finite dimensional representation of $\g$ is one dimensional.
Furthermore, the set of isomorphism classes of irreducible finite dimensional representations of $\g$ is parametrized by
$(\g/\C^{1}(\g))^{*}$.
\end{lema}
\noindent\textbf{Proof.}
Let $V$ be an irreducible finite dimensional $\g$-module of dimension $n \geq 1$ defined by $\phi : \g \rightarrow \gl(V)$. 
Since $V$ is simple, Lie's theorem tells us that $V$ contains a common eigenvector $v$ for all the endomorphisms in $\g$ 
(cf. \cite{Hum1}, Thm. 4.1). 
Then, the non trivial submodule $k.v \subset V$ should coincide with $V$, for $V$ is irreducible, so $\dim_{k}(V) = 1$.  

On the other hand, since $\gl(V) \simeq k$, the morphism $\phi : \g \rightarrow \gl(V)$ takes values in an abelian Lie algebra, 
so $\C^{1}(\g) \subset \Ker(\phi)$, which in turn implies that $\phi$ induces a morphism $\bar{\phi} : \g/\C^{1}(\g) \rightarrow k$. 
Conversely, given a linear form $\bar{\phi} \in (\g/\C^{1}(\g))^{*}$, we have a Lie algebra homomorphism 
$\phi : \g \rightarrow k \simeq \gl(k)$, whose kernel contains $\C^{1}(\g)$. 
Thus, $k$ turns out to be a irreducible representation of $\g$.
The lemma is proved.
\qed

From the previous lemma and Proposition \ref{prop:mod}, and taking into account that 
an irreducible solvable finite dimensional $\ym(n)$-module is also nilpotent, we may obtain the following theorem:
\begin{teo}
Every irreducible solvable finite dimensional representation of the Yang-Mills algebra $\ym(n)$ is of dimension $1$.
Moreover, the set of isomorphism classes of irreducible solvable finite dimensional representations of the Yang-Mills algebra $\ym(n)$ is
parametrized by $(\ym(n)/C^{1}(\ym(n)))^{*}$. 
\end{teo}
%

\begin{rem}
The previous theorem is not only true for Yang-Mills algebras but also for any Lie algebra $\g$
provided with a locally finite dimensional $\NN$-grading, taking into account that the set of isomorphism classes of irreducible
solvable finite dimensional representations is parametrized by $(\g/C^{1}(\g))^{*}$.
\end{rem}

\section{\texorpdfstring {The ideal $\tym(n)$}
        {The ideal tym(n)}}
\label{sec:tym}

In this section we shall study in detail the ideal $\tym(n)$ of the Lie Yang-Mills algebra and its associative version $\TYM(n)$. 
This ideal is the key point of the construction of the family of representations we shall define. 
In order to achieve this construction we need to prove that $\TYM(n)$ is a free algebra. 
We shall then perform detailed computations in order to prove this fact, for which we shall make use of the bar construction for
augmented differential graded algebras (or A$_{\infty}$-algebras). 
We suggest \cite{Lef1, Kel2} as a reference. 
Although some of the results of this section are mentioned in \cite{MS1}, our aim here is to give detailed proofs of the results that we will need later. 

The starting point is to define $k$-linear morphisms $d_{i}$, $i = 1, \dots, n$, given by
\begin{equation}
\label{eq:derivacionesdev}
\begin{split}
   d_{i} : V(n) &\rightarrow V(n)
   \\
   d_{i} (x_{j}) &= \delta_{i j}.
\end{split}
\end{equation}
They can be uniquely extended to derivations $d_{i}$, $i= 1, \dots, n$, on $TV(n)$.
Since
\begin{equation}
\label{eq:derivacionesYangMills}
     d_{i}([x_{j},x_{k}]) = d_{i}(x_{j} x_{k} - x_{k} x_{j})
                            = \delta_{i j} x_{k} + x_{j} \delta_{i k} - \delta_{i k} x_{j} - x_{k} \delta_{i j} = 0, \forall \hskip 0.6mm
                            i, j, k = 1, \dots, n,
\end{equation}
we see that each $d_{i}$ induces a derivation on $TV(n)/\cl{R(n)} = \YM(n)$, which we will also denote $d_{i}$.

The following proposition characterizes the algebra $\U(\tym(n))$ as a subalgebra of $\U(\ym(n)) = \YM(n)$.
\begin{prop}
\label{prop:nucleotym}
The inclusion $\mathrm{inc} : \tym(n) \hookrightarrow \ym(n)$ induces a monomorphism of algebras
$\U(\mathrm{inc}) : \U(\tym(n)) \hookrightarrow \U(\ym(n))$, with image $\bigcap_{i=1}^{n} \Ker(d_{i})$.
\end{prop}
\noindent\textbf{Proof.}
The first statement is a direct consequence of the Poincar\'e-Birkhoff-Witt Theorem (cf. \cite{Dix1}, Sec. 2.2.6, Prop. 2.2.7).
As it is usual, we will identify $\U(\tym(n))$ with its image by $\U(\mathrm{inc})$ in $\U(\ym(n))$.

Let us prove the second statement.
On the one hand, if $z \in \tym(n) = [\ym(n),\ym(n)]$, then $d_{i} (z) = 0$, $i = 1, \dots, n$.
Since $\tym(n)$ generates the algebra $\U(\tym(n))$ and $d_{i}$ is a derivation, each $d_{i}$ vanishes in $\U(\tym(n))$.
Hence,
\[     \U(\tym(n)) \subseteq \bigcap_{i=1}^{n} \Ker(d_{i}).     \]

We next choose an ordered basis of $\tym(n)$ as $k$-vector space, denoted by $\B' = \{ y_{j} : j \in J \}$.
As a consequence, the set $\B = \{ x_{1}, \dots, x_{n} \} \cup \B'$ is an ordered basis of $\ym(n)$.
By the Poincar\'e-Birkhoff-Witt Theorem, given $z \in \U(\ym(n))$,
\[     z = \underset{\text{\begin{tiny}$\begin{matrix}
                 j_{1},\dots,j_{l} \in J \hskip 0.5mm \text{not equal}
                 \\
                 (r_{1}, \dots,r_{n},s_{1},\dots, s_{l}) \in \NN_{0}^{n+l}
                 \end{matrix}$\end{tiny}}}{\sum}
                 c_{(r_{1}, \dots,r_{n},s_{1},\dots, s_{l})}^{j_{1},\dots,j_{l}}
           x_{1}^{r_{1}} \dots x_{n}^{r_{n}} y_{j_{1}}^{s_{1}} \dots y_{j_{l}}^{s_{l}}.     \]
Since $d_{i} (x_{j}^{p}) = p x_{j}^{p-1} \delta_{i j}$ and $d_{i} (y_{j_{m}}) = 0$, $m = 1, \dots, l$,
\[     d_{i}(z) = \underset{\text{\begin{tiny}$\begin{matrix}
                 j_{1},\dots,j_{l} \in J \hskip 0.5mm \text{not equal}
                 \\
                 (r_{1}, \dots,r_{n},s_{1},\dots, s_{l}) \in \NN_{0}^{n+l}
                 \end{matrix}$\end{tiny}}}{\sum}
                 c_{(r_{1}, \dots,r_{n},s_{1},\dots, s_{l})}^{j_{1},\dots,j_{l}}
                 r_{i} x_{1}^{r_{1}} \dots x_{i}^{r_{i}-1} \dots x_{n}^{r_{n}} y_{j_{1}}^{s_{1}} \dots y_{j_{l}}^{s_{l}}.     \]

Let us suppose that $d_{i}(z) = 0$, for all $i = 1, \dots, n$.
By the Poincar\'e-Birkhoff-Witt Theorem, we obtain that $r_{i} = 0$, for all $i = 1, \dots, n$. 
This in turn implies that $z \in \U(\tym(n))$, and the other inclusion is proved. 
\qed

We shall consider the following filtration on the algebra $\YM(n)$
\[     F^{j} = \begin{cases}
                  0 &\text{if $j = 0$,}
                  \\
                  \{ z \in \TYM(n) : d_{i} (z) \in F^{j-1}, \forall \hskip 0.5mm i, 1 \leq i \leq n \} &\text{if $j \in \NN$.}
               \end{cases}
\]
By the previous proposition, $F^{1} = \TYM(n)$.

\begin{lema}
The filtration $\{ F^{j} \}_{j\in \NN_{0}}$ defined on $YM(n)$ is increasing, multiplicative, exhaustive, Hausdorff and such that
$x_{i} \in F^{2}$, for all $i = 1, \dots, n$.
\end{lema}
\noindent\textbf{Proof.}
Cf. \cite{MS1}, Lemma 28.
\qed

The following result is implicit in the analysis of \cite{MS1}.
\begin{lema}
\label{lema:potencial}
Let $A$ be an associative $k$-algebra with unit.
Let $\partial_{i}$, $i = 1, \dots, n$ be the usual derivations on $k[t_{1}, \dots, t_{n}]$ and define the derivations
$D_{i} = \id_{A} \otimes \partial_{i}$ on the algebra $A [t_{1}, \dots, t_{n}] \simeq A \otimes_{k} k[t_{1}, \dots, t_{n}]$.
Then, given polynomials $p_{1}, \dots, p_{n} \in A[t_{1}, \dots, t_{n}]$ such that $D_{i}(p_{j}) = D_{j}(p_{i})$, for all $i, j = 1, \dots, n$,
there exists a polynomial $P \in A[t_{1}, \dots, t_{n}]$ such that $D_{i}(P) = p_{i}$, $i = 1, \dots, n$.
\end{lema}
\noindent\textbf{Proof.}
The classical proof for $A = k$ (cf. \cite{Cou1}, Lemma 2.2) works as well in this context.
\qed

Let $\mathrm{Gr}_{F^{\bullet}}(\YM(n))$ be the associated graded algebra of $\YM(n)$ provided with the filtration $\{ F^{j} \}_{j\in \NN_{0}}$,
and let us denote $\bar{z}$ the class in $\mathrm{Gr}_{F^{\bullet}}(\YM(n))$ of an element $z \in \YM(n)$.

%

The derivations $d_{i}$, $i = 1, \dots, n$, induce morphisms on the associated graded algebra,
which are also derivations and we shall denote them in the same way.

Furthermore, $F^{1}/F^{0}$ is a subalgebra of $\mathrm{Gr}_{F^{\bullet}}(\YM(n))$, isomorphic as an algebra to $F^{1} = \TYM(n)$.
From now on, we shall make use of this identification.

\begin{lema}
\label{lema:tecnicotym}
The algebra $\mathrm{Gr}_{F^{\bullet}}(\YM(n))$ satisfies the following properties:
\begin{itemize}
\item[(i)] The elements $\bar{x}_{i}$, $i = 1, \dots, n$, commute with each other and with the subalgebra $F^{1} / F^{0}$.

\item[(ii)] The elements $\bar{x}_{i}$, $i = 1, \dots, n$, and $F^{1} / F^{0}$ generate $\mathrm{Gr}_{F^{\bullet}}(\YM(n))$.

\item[(iii)] The algebra generated by the elements $\bar{x}_{i}$ and $F^{1} / F^{0}$ is isomorphic to
$(F^{1} / F^{0}) \otimes_{k} k[t_{1}, \dots, t_{n}]$.
Hence, there exists an isomorphism of algebras $(F^{1} / F^{0}) \otimes_{k} k[t_{1}, \dots, t_{n}] \simeq \mathrm{Gr}_{F^{\bullet}}(\YM(n))$.
\end{itemize}
\end{lema}
\noindent\textbf{Proof.}
Cf. \cite{MS1}, Lemma 29.
Observe that Lemma \ref{lema:potencial} is necessary to prove (ii) and (iii).
\qed

The projection $\ym(n) \rightarrow \ym(n)/\tym(n) \simeq V(n)$ provides an action of the Lie Yang-Mills algebra $\ym(n)$ on the symmetric algebra
$S(V(n))$ such that the restricted action of $\tym(n)$ on $S(V(n))$ is trivial.

In the same way as we have done for the Yang-Mills algebra, we define the \emph{special grading} of $S(V(n))$:
we consider $V(n)$ concentrated in degree $2$ and we identify $S(V(n))$ and $S_{gr}(V(n))$, where the latter is the symmetric algebra in the
category of graded $k$-vector spaces.
The \emph{usual grading} of $S(V(n))$ is given by considering $V(n)$ concentrated in degree $1$.

If the Yang-Mills algebra $\ym(n)$ and the symmetric algebra $S(V(n))$ are provided with the special grading,
then $S(V(n))$ is a graded module over $\ym(n)$.

We will compute the cohomology of $\ym(n)$ with coefficients on $S(V(n))$ because it will be used in order to prove weak convergence of the
spectral sequence defined in Corollary \ref{coro:sucesionespectral}.
\begin{prop}
\label{prop:homologiayangmillssimetrica}
The Lie homology of $\ym(n)$ with coefficients in the module $S(V(n))$ is given by
\[
H_{\bullet} (\ym(n),S(V(n))) = \begin{cases}
                               k, &\text{if $\bullet = 0$,}
                               \\
                               0, &\text{if $\bullet \geq 2$.}
                               \end{cases}
\]
The homology in degree $1$ is the direct sum of vector spaces $H_{1}^{p}$ ($p \in \NN_{0}$), where
\[     \dim_{k} (H_{1}^{p}) = \begin{cases}
                              n \frac{(n+p-1)!}{(n-1)!} - \frac{(n+p)!}{(n-1)! (p+1)!}, &\text{if $p = 0$ or $p = 1$.}
                              \\
                              n \frac{(n+1)!}{(n-1)! 2!} - \frac{(n+2)!}{(n-1)! 3!} - n, &\text{if $p = 2$,}
                              \\
                              n \frac{(n+p-1)!}{(n-1)! p!} - \frac{(n+p)!}{(n-1)! (p+1)!} - n \frac{(n+p-3)!}{(n-1)! (p-2)!}
                              + \frac{(n+p-4)!}{(n-1)! (p-3)!}, &\text{if $p \geq 3$,}
                              \end{cases}. \]
\end{prop}
\noindent\textbf{Proof.}
We shall use the Koszul resolution (9) of $\U(\ym(n))$ described by \cite{CD1} for homological computations.
The Koszul complex $(C_{\bullet} (\YM(n),S(V(n))),d_{\bullet})$ 
computing the homology $\Tor^{\U(\ym(n))}_{\bullet} (k,S(V(n))) \simeq H_{\bullet} (\ym(n), S(V(n)))$
is given by tensoring the Koszul resolution (9) of \cite{CD1} with $S(V(n))$ over $\YM(n)$ 
\begin{equation}
\label{eq:complejocohomologiayangmillssimetrica}
    0 \longrightarrow S(V(n))[-4] \overset{d_{3}}{\longrightarrow} (S(V(n)) \otimes V(n))[-2] \overset{d_{2}}{\longrightarrow} S(V(n)) \otimes V(n)
    \overset{d_{1}}{\longrightarrow} S(V(n)) \longrightarrow 0,
\end{equation}
where the differentials are
\begin{align*}
   d_{3}(w) &= \sum_{i=1}^{n} x_{i}.w \otimes x_{i},
   \\
   d_{2}(w \otimes x_{i}) &= \sum_{j=1}^{n} (x_{j}^{2}.w \otimes x_{i} - x_{i} x_{j}.w \otimes x_{j}),
   \\
   d_{1}(w \otimes x_{i}) &= x_{i}.w,
\end{align*}
for $S(V(n))$ is commutative.
We have shifted some terms of the complex so the complex \eqref{eq:complejocohomologiayangmillssimetrica} is composed of homogeneous morphisms of degree $0$.

We immediately get that $H_{\bullet} (\ym(n), S(V(n))) = 0$ if $\bullet > 3$.

The complex \ref{eq:complejocohomologiayangmillssimetrica} is the direct sum of the following subcomplexes 
of finite dimensional $k$-vector spaces
\begin{equation}
    0 \longrightarrow S^{p-1}(V(n))[-4] \overset{d_{3}^{p-1}}{\longrightarrow} (S^{p}(V(n)) \otimes V(n))[-2] \overset{d_{2}^{p}}{\longrightarrow}
    S^{p+2}(V(n)) \otimes V(n) \overset{d_{1}^{p+2}}{\longrightarrow} S^{p+3}(V(n)) \longrightarrow 0,
\end{equation}
where $p \in \ZZ$ and we consider $S^{p}(V(n)) = 0$, if $p < 0$.
We define
\[     H_{\bullet}^{p} (\ym(n),S(V(n))) = \begin{cases}
                                          \Ker (d_{\bullet}^{p})/\mathrm{Im}(d_{\bullet + 1}^{p-1}), &\text{if $\bullet \neq 1$,}
                                          \\
                                          \Ker (d_{\bullet}^{p})/\mathrm{Im}(d_{\bullet + 1}^{p-2}), &\text{if $\bullet = 1$.}
                                          \end{cases}
\]

Notice that $d_{3}^{p}$ is injective for $p \in \NN_{0}$. 
This is proved as follows. 
If $d_{3}^{p} (w) = 0$, then $x_{i} . w = 0$, for all $i = 1, \dots, n$, so $w = 0$, 
because $\{ x_{i} \}_{i = 1, \dots, n}$ is a basis of $V(n)$ and $S(V(n))$ is entire.
Hence, $H_{3} (\ym(n), S(V(n))) = 0$.

On the other hand, the morphism $d_{1}^{p}$ is surjective for $p \in \NN_{0}$, since, if $w \in S^{p+1}(V(n))$, then
there exists $i$ such that $w = x_{i}.w'$ ($p +1 > 0$!), i.e. $w = d_{1} (w')$.
The morphism $d_{1}$ being homogeneous of degree $1$, if $w \in S^{0}(V(n)) = k$, there does not exist $w' \in S(V(n)) \otimes V(n)$
such that $d_{1}(w') = w$.
As a consequence, $H_{0} (\ym(n), S(V(n))) = k$.

We will now show that $H_{2} (\ym(n), S(V(n))) = 0$.
Let now $w = \sum_{i=1}^{n} w_{i} \otimes x_{i} \in S^{p}(V(n)) \otimes V(n)$ (where $w_{i} \in S^{p}(V(n))$ for $i = 1, \dots, n$)
be in the kernel of the morphism $d_{2}^{p}$.
Then $ 0 = d_{2}(w) = \sum_{i,j=1}^{n} (w_{i} x_{j}^{2} \otimes x_{i} - w_{i} x_{i} x_{j} \otimes x_{j})$, and using that
$\{ x_{i} \}_{i = 1, \dots, n}$ is a basis of $V(n)$, it turns out that
$\sum_{j=1}^{n} (w_{i} x_{j}^{2} - w_{j} x_{i} x_{j}) = 0$, for $i = 1, \dots, n$.
This is equivalent to
\[     w_{i} \sum_{j=1}^{n} x_{j}^{2} = x_{i} \sum_{j=1}^{n} x_{j} w_{j}, \hskip 0.5mm i = 1, \dots, n.     \]
Since $S(V(n))$ is a unique factorization domain and the elements $x_{i}$ are prime, this identity implies that $x_{i}$ divides $w_{i}$,
for $i = 1, \dots, n$.

Let $w'_{i}$ be such that $w_{i} = x_{i} w'_{i}$.
We can rewrite the previous equation as follows
\[     \sum_{j=1}^{n} x_{j}^{2} (w'_{i} - w'_{j}) = 0, \forall \hskip 0.5mm i = 1, \dots, n.     \]
Fix $i_{1}, i_{2}$ such that $1 \leq i_{1} < i_{2} \leq n$. 
Then, 
$\sum_{j = 1}^{n} x_{j}^{2} (w'_{i_{1}} - w'_{i_{2}}) = 0$, and $S(V(n))$ being entire, we see that
$w'_{i_{1}} = w'_{i_{2}}$, for all $i_{1}, i_{2}, 1 \leq i_{1} < i_{2} \leq n$.
Let us call this element $w'$.
Hence, $w_{i} = x_{i} w'$.

Now, $d_{3} (w') = \sum_{i=1}^{n} w' x_{i} \otimes x_{i} = \sum_{i=1}^{n} w_{i} \otimes x_{i} = w$.
We conclude that $H_{2} (\ym(n), S(V(n))) = 0$.

We finally compute the homology in degree $1$.

Since $d_{1}^{p}$ is surjective, we have that
\[     \dim_{k} (\Ker(d_{1}^{p})) = \dim_{k} (S^{p}(V(n)) \otimes V(n)) - \dim_{k} (S^{p+1}(V(n)))
                                 = n \frac{(n+p-1)!}{(n-1)! p!} - \frac{(n+p)!}{(n-1)! (p+1)!}.     \]
On the other hand, as $H_{2}(\ym(n), S(V(n))) = 0$ we know that $\Ker (d_{2}^{p-2}) = \mathrm{Im} (d_{3}^{p-3})$.
Moreover, injectivity of $d_{3}^{p-3}$ yields that
\[     \dim_{k} (\mathrm{Im}(d_{2}^{p-2})) = \dim_{k} (S^{p-2}(V(n)) \otimes V(n)) - \dim_{k} (S^{p-3}(V(n)))
                                 = n \frac{(n+p-3)!}{(n-1)! (p-2)!} - \frac{(n+p-4)!}{(n-1)! (p-3)!},     \]
for $p \geq 3$.

If $p = 0, 1$, then $\mathrm{Im} (d_{2}^{p-2}) = \{ 0 \}$ and in case $p = 2$, $S^{p-3}(V(n)) = \{ 0 \}$, so
\[     \dim_{k} (\mathrm{Im}(d_{2}^{p-2})) = \dim_{k} (S^{p-2}(V(n)) \otimes V(n)) = n.     \]

The proposition is proved.
\qed

\begin{coro}
\label{coro:sucesionespectral}
The filtration $\{ F^{p}C_{\bullet} (\YM(n),S(V(n))) \}_{p \in \ZZ}$ of the complex $(C_{\bullet} (\YM(n),S(V(n))),d_{\bullet})$ of 
\eqref{eq:complejocohomologiayangmillssimetrica} given by
\begin{multline}
\label{eq:complejocohomologiayangmillssimetricafiltracion}
    F^{p}C_{\bullet} (\YM(n),S(V(n)))
    \\
    = (0 \longrightarrow S^{\geq -p}(V(n)) \overset{d_{3}}{\longrightarrow} S^{\geq -p}(V(n)) \otimes V(n)
    \overset{d_{2}}{\longrightarrow} S^{\geq -p}(V(n)) \otimes V(n) \overset{d_{1}}{\longrightarrow} S^{\geq -p}(V(n)) \longrightarrow 0),
\end{multline}
is increasing, exhaustive and Hausdorff, and the spectral sequence associated to this filtration weakly converges to the homology
$H_{\bullet}(\ym(n) , S(V(n)))$ of the complex $(C_{\bullet} (\YM(n),S(V(n))),d_{\bullet})$.
\end{coro}
\noindent\textbf{Proof.}
We shall indicate the consecutive steps of the spectral sequence, following the standard construction detailed in \cite{Wei1}, Sec. 5.4.
First, notice that, since $d_{1}$ and $d_{3}$ are homogeneous morphisms of degree $1$, and $d_{2}$ is homogeneous of degree $2$,
then
\[     d_{i} (F^{p}C_{\bullet} (\YM(n),S(V(n)))) \subseteq F^{p-1}C_{\bullet} (\YM(n),S(V(n))),     \]
for $i = 1, 2, 3$, and
\[     d_{2} (F^{p}C_{\bullet} (\YM(n),S(V(n)))) \subseteq F^{p-2}C_{\bullet} (\YM(n),S(V(n))).     \]
This in turns implies that the differentials $d^{0}_{p,q}$ are $0$.
Besides, as
\[     E^{0}_{p,q} = F^{p} C_{p+q}(\ym(n),S(V(n)))/F^{p-1}C_{p+q}(\ym(n),S(V(n))),     \]
the spectral sequence is concentrated in the set of $(p,q)$ such that $0 \leq p + q \leq 3$, and $p \leq 0$, and
\[     E^{0}_{p,q} = \begin{cases}
                     S^{-p}(V(n)), &\text{if $q = -p$,}
                     \\
                     S^{-p}(V(n)) \otimes V(n), &\text{if $q = -p+1$,}
                     \\
                     S^{-p}(V(n)) \otimes V(n), &\text{if $q = -p+2$,}
                     \\
                     S^{-p}(V(n)), &\text{if $q = -p+3$,}
                     \\
                     0 &\text{in other case.}
                     \end{cases}
\]
We may picture the terms $E^{0}_{\bullet, \bullet}$ as follows
\begin{center}
\begin{figure}[H]
\[
\xymatrix@R-12pt
{
&
&
&
&
&
\\
\bullet
\ar[d]^{0}
&
\bullet
\ar[d]^{0}
&
\bullet
\ar[d]^{0}
&
0
&
0
&
E_{\bullet,\bullet}^{0}
\\
\bullet
\ar[d]^{0}
&
\bullet
\ar[d]^{0}
&
\bullet
\ar[d]^{0}
&
\bullet
\ar[d]^{0}
\ar@{.}[ul]
&
0
&
\\
0
&
\bullet
\ar[d]^{0}
&
\bullet
\ar[d]^{0}
&
\bullet
\ar[d]^{0}
&
0
&
\\
0
&
0
&
\bullet
\ar[d]^{0}
&
\bullet
\ar[d]^{0}
&
0
&
\\
0
&
0
&
0
&
\bullet
\ar[uuuuu]^>{q}
\ar[rr]^>{p}
\ar@{.}[uuulll]
&
0
&
}
\]
\caption[margin=150pt]{Zeroth term $E^{0}_{\bullet,\bullet}$ of the spectral sequence.
The dotted lines indicate the limits wherein the spectral sequence is concentrated.}
\end{figure}
\end{center}

Since $d^{0}_{p,q} = 0$ and $E^{0}_{p,q} = E^{1}_{p,q}$
, it turns out that
\[     d^{1}_{p,q} = \begin{cases}
                     0, &\text{if $q = -p$,}
                     \\
                     d_{1}^{-p}, &\text{if $q = -p+1$,}
                     \\
                     0, &\text{if $q = -p+2$,}
                     \\
                     d_{3}^{-p}, &\text{if $q = -p+3$,}
                     \\
                     0 &\text{in other case.}
                     \end{cases}
\]

\begin{center}
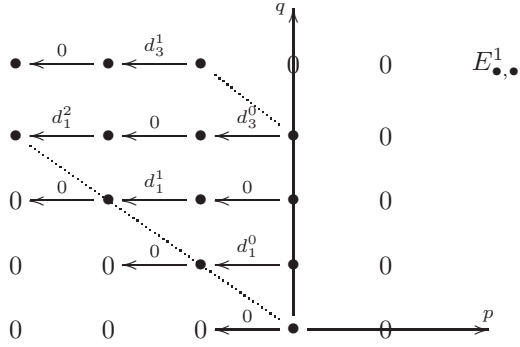
\begin{figure}[H]
\[
\xymatrix@R-12pt
{
&
&
&
&
&
\\
\bullet
&
\bullet
\ar[l]_{0}
&
\bullet
\ar[l]_{d_{3}^{1}}
&
0
&
0
&
E_{\bullet,\bullet}^{1}
\\
\bullet
&
\bullet
\ar[l]_{d_{1}^{2}}
&
\bullet
\ar[l]_{0}
&
\bullet
\ar[l]_{d_{3}^{0}}
\ar@{.}[ul]
&
0
&
\\
0
&
\bullet
\ar[l]_{0}
&
\bullet
\ar[l]_{d_{1}^{1}}
&
\bullet
\ar[l]_{0}
&
0
&
\\
0
&
0
&
\bullet
\ar[l]_{0}
&
\bullet
\ar[l]_{d_{1}^{0}}
&
0
&
\\
0
&
0
&
0
&
\bullet
\ar[uuuuu]^>{q}
\ar[rr]^>{p}
\ar@{.}[uuulll]
\ar[l]_{0}
&
0
&
}
\]
\caption[margin=150pt]{First term $E^{1}_{\bullet,\bullet}$ of the spectral sequence.
}
\end{figure}
\end{center}

In consequence, the second step of the spectral sequence is
\[     E^{2}_{p,q} = \begin{cases}
                     S^{-p}(V(n))/\mathrm{Im}(d_{1}^{-p+1}), &\text{if $q = -p$,}
                     \\
                     \Ker(d_{1}^{-p}), &\text{if $q = -p+1$,}
                     \\
                     \mathrm{Im}(d_{3}^{-p+1}), &\text{if $q = -p+2$,}
                     \\
                     \Ker(d_{3}^{-p}), &\text{if $q = -p+3$,}
                     \\
                     0 &\text{in other case,}
                     \end{cases}
\]
and the differentials are 
\[     d^{2}_{p,q} = \begin{cases}
                     d_{2}^{-p}, &\text{if $q = -p+2$,}
                     \\
                     0 &\text{if not.}
                     \end{cases}
\]
In this case,
\begin{center}
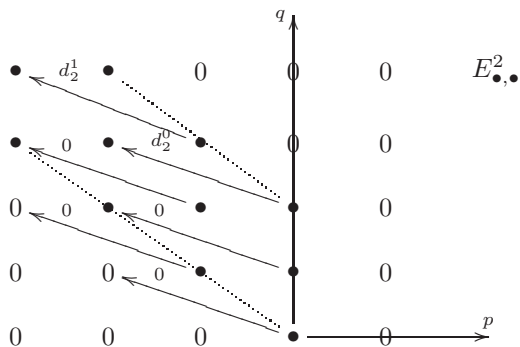
\begin{figure}[H]
\[
\xymatrix@R-12pt
{
&
&
&
&
&
\\
\bullet
&
\bullet
&
0
&
0
&
0
&
E^{2}_{\bullet,\bullet}
\\
\bullet
&
\bullet
&
\bullet
\ar[llu]_(.75){d_{2}^{1}}
&
0
&
0
&
\\
0
&
\bullet
&
\bullet
\ar[llu]_(.75){0}
&
\bullet
\ar[llu]_(.75){d_{2}^{0}}
\ar@{.}[uull]
&
0
&
\\
0
&
0
&
\bullet
\ar[llu]_(.75){0}
&
\bullet
\ar[llu]_(.75){0}
&
0
&
\\
0
&
0
&
0
&
\bullet
\ar[llu]_(.75){0}
\ar[uuuuu]^>{q}
\ar[rr]^>{p}
\ar@{.}[uuulll]
&
0
&
}
\]
\caption[margin=150pt]{Step $E^{2}_{\bullet, \bullet}$ of the spectral sequence.
}
\end{figure}
\end{center}
As a consequence, the third step of the spectral sequence is
\[     E^{3}_{p,q} = \begin{cases}
                     H_{0}^{-p}(\ym(n),S(V(n))), &\text{if $q = -p$,}
                     \\
                     H_{1}^{-p}(\ym(n),S(V(n))), &\text{if $q = -p+1$,}
                     \\
                     H_{2}^{-p}(\ym(n),S(V(n))), &\text{if $q = -p+2$,}
                     \\
                     H_{3}^{-p}(\ym(n),S(V(n))), &\text{if $q = -p+3$,}
                     \\
                     0 &\text{in other case.}
                     \end{cases}
\]
Hence, the spectral sequence is weakly convergent because of Proposition \ref{prop:homologiayangmillssimetrica}.
\qed

We shall consider the graded algebra $\YM(n) \otimes \Lambda^{\bullet} V(n)$ with the usual multiplication and the grading given by
putting $\YM(n)$ in degree zero and the usual grading of $\Lambda^{\bullet} V(n)$.
We define a differential $d$ of degree $1$ on this algebra by the formula
\[     d (z \otimes w) = \sum_{i=1}^{n} d_{i}(z) \otimes (x_{i} \wedge w),      \]
where $z \in \YM(n)$ and $w \in \Lambda^{\bullet} V(n)$.
The identity $d \circ d = 0$ is immediate. 
The map $d$ is a graded derivation, since
\begin{align*}
   d((z \otimes w) (z' \otimes w')) &= d(z z' \otimes w \wedge w') = \sum_{i=1}^{n} d_{i}(z z') \otimes (x_{i} \wedge w \wedge w')
                                        = \sum_{i=1}^{n} (d_{i}(z) z' + z d_{i}(z')) \otimes (x_{i} \wedge w \wedge w')
   \\
                                   &= (\sum_{i=1}^{n} d_{i}(z) \otimes (x_{i} \wedge w)) (z' \otimes w') + (-1)^{|w|}
                                      (z \otimes w)(\sum_{i=1}^{n} d_{i}(z')) \otimes (x_{i} \wedge w')
   \\
                                   &= d(z \otimes w) (z' \otimes w') + (-1)^{|z \otimes w|}(z \otimes w) d(z' \otimes w').
\end{align*}

If $\epsilon_{\ym(n)}$ is the augmentation of $\YM(n)$ and $\epsilon'$ is the augmentation of the exterior algebra $\Lambda^{\bullet} V(n)$,
we shall consider the augmentation of $\YM(n) \otimes \Lambda^{\bullet} V(n)$ given by the usual formula
$\epsilon = \epsilon_{\ym(n)} \otimes \epsilon'$.
Hence, we have defined a structure of augmented differential graded algebra on $\YM(n) \otimes \Lambda^{\bullet} V(n)$.

The next result says that this augmented differential graded algebra is in fact a model for $\TYM(n)$.
\begin{prop}
\label{prop:quasiiso}
If we consider the algebra $\TYM(n)$ as an augmented differential graded algebra concentrated in degree zero,
with zero differential and augmentation $\epsilon_{\tym(n)}$, then the morphism
\begin{align*}
   \mathrm{inc} : \TYM(n) &\rightarrow \YM(n) \otimes \Lambda^{\bullet} V(n)
   \\
   z &\mapsto z \otimes 1
\end{align*}
is a quasi-isomorphism of augmented differential graded algebras and we will write $\TYM(n) \simeq_{q} \YM(n) \otimes \Lambda^{\bullet} V(n)$.
\end{prop}
\noindent\textbf{Proof.}
The map $\mathrm{inc}$ is a morphism of graded algebras, as we can easily see.
It also commutes with differentials by Proposition \ref{prop:nucleotym}.
Furthermore, the same proposition also implies that $\mathrm{inc}$ induces an isomorphism between $\TYM(n)$ and
$H^{0} (\YM(n) \otimes \Lambda^{\bullet} V(n)) = Z^{0} (\YM(n) \otimes \Lambda^{\bullet} V(n))$,
since $z \in Z^{0} (\YM(n) \otimes \Lambda^{\bullet} V(n))$ if and only if
$z = v \otimes 1$, with $v \in \YM(n)$ and $d(z) = \sum_{i=1}^{n} d_{i}(v) \otimes x_{i} = 0$,
which in turn happens if and only if $d_{i}(v) = 0$, for all $i = 1, \dots, n$.

Note that $\mathrm{inc}$ also commutes with the augmentations.

We shall now proceed to prove that $\mathrm{inc}$ induces an isomorphism in cohomology.
It is thus necessary to compute the cohomology of the underlying cochain complex $(\YM(n) \otimes \Lambda^{\bullet} V(n),d)$.
We will write the cochain complex $(\YM(n) \otimes \Lambda^{\bullet} V(n),d)$ as a chain complex
in the usual way $C_{\bullet} = \YM(n) \otimes \Lambda^{n - \bullet} V(n)$.

We consider $(C_{\bullet},d)$ provided with the filtration $\{ F_{\bullet}C \}_{\bullet \in \ZZ}$ defined as follows
\[     F_{p}C_{q} = F^{p+n-q} \otimes \Lambda^{n-q} V(n).     \]
Notice that $\{ F_{\bullet}C \}_{\bullet \in \ZZ}$ is an increasing, bounded below and exhaustive filtration and
$d(F_{p}C_{q}) \subseteq F_{p-2}C_{q-1}$.
Hence, $\{ F_{\bullet}C \}_{\bullet \in \ZZ}$ is a filtration of complexes and in turn it induces a spectral sequence whose second term is 
\begin{align*}
  E^{2}_{p,q} = F_{p} C_{p+q}/F_{p-1} C_{p+q}
  &= (F^{n-q} \otimes \Lambda^{n-(p+q)} V(n))/ (F^{n-q-1} \otimes \Lambda^{n-(p+q)} V(n))
  \\
  &\simeq (F^{n-q}/F^{n-q-1}) \otimes \Lambda^{n-(p+q)} V(n)
  = \mathrm{Gr}_{F^{\bullet}}(YM(n))^{n-q} \otimes \Lambda^{n-(p+q)} V(n)
  \\
  &\simeq \TYM(n) \otimes S^{n-q}(V(n)) \otimes \Lambda^{n-(p+q)} V(n),
\end{align*}
where the last isomorphism follows from the last item of Lemma \ref{lema:tecnicotym}.
The differential $d^{2}_{p,q} : E^{2}_{p,q} \rightarrow E^{2}_{p-2,q+1}$ can be written in the following simple way
\[
\xymatrix@C60pt
{
 E^{2}_{p,q}
 \ar^{d^{2}_{p,q}}[r] 
 \ar@{=}[d]
 &
 E^{2}_{p-2,q+1}
 \ar@{=}[d]
 \\
 (F^{n-q}/F^{n-q-1}) \otimes \Lambda^{n-(p+q)} V(n)
 \ar^{\bar{d}}[r]
 \ar@{=}[d]
 &
 (F^{n-q-1}/F^{n-q-2}) \otimes \Lambda^{n-(p+q)+1} V(n)
 \ar@{=}[d]
 \\
 \TYM(n) \otimes S^{n-q}(V(n)) \otimes \Lambda^{n-(p+q)} V(n)
 \ar^{d'_{p,q}}[r]
 &
 \TYM(n) \otimes S^{n-q}(V(n)) \otimes \Lambda^{n-(p+q)+1} V(n)
}
\]
where $\bar{d}$ is the morphism induced by $d$ and $d'_{p,q}$ is the map given by
\[     d'_{p,q}(z \otimes v \otimes w) = \sum_{i=1}^{n} z \otimes \partial_{i}(v) \otimes (x_{i} \wedge w).     \]

We see that this complex is exact except in case $p = 0$ and $q = n$.
This follows from the fact that the differential $d'_{p,q}$ is the $\TYM(n)$-linear extension of the differential of the de Rham complex of the
algebra $S(V(n))$, whose cohomology is zero, except in degree zero, where it equals $k$ (cf. \cite{Wei1}, Coro. 9.9.3).
This implies that the spectral sequence collapses in the third step, for the unique non zero element is $E^{3}_{0,n} = \TYM(n)$.
As the filtration is bounded below and exhaustive, the Classical Convergence Theorem tells us that the spectral sequence is convergent 
and it converges to the homology of the complex $(C_{\bullet},d_{\bullet})$ (cf. \cite{Wei1}, Thm. 5.5.1).

Finally, $H^{\bullet} (\YM(n) \otimes \Lambda^{\bullet} V(n)) = H_{n-\bullet}(C) = 0$, if $\bullet \neq 0$,
and $H^{0} (\YM(n) \otimes \Lambda^{\bullet} V(n)) = H_{n}(C) = \TYM(n)$.
\qed

Let $B^{+}(\place)$ denote the bar construction for augmented differential graded algebras or more generally, for A$_{\infty}$-algebras
(cf. \cite{Lef1}, Notation 2.2.1.4).
From the previous proposition we obtain that $B^{+}(\TYM(n)) \simeq_{q} B^{+}(\YM(n) \otimes \Lambda^{\bullet} V(n))$.
Since the underlying complex of the bar complex of an augmented (graded) algebra coincides with the normalized (graded) Hochschild bar complex
with coefficients in the (graded) bimodule $k$ (cf. \cite{Lef1}, Proof of Lemma 2.2.1.9),
the respective Hochschild homologies of differential graded algebras with coefficients in $k$ are isomorphic, i.e.
\[     H_{\bullet} (\TYM(n) , k) = H_{\bullet} (B^{+}(\TYM(n)))
       \simeq H_{\bullet} (B^{+}(\YM(n) \otimes \Lambda^{\bullet} V(n)))
       = H_{\bullet} (\YM(n) \otimes \Lambda^{\bullet} V(n),k).     \]

We define another filtration on the augmented differential graded algebra $\YM(n) \otimes \Lambda V(n)$ by the formula
\begin{equation}
\label{eq:anotherfiltration}
     F_{p} (\YM(n) \otimes \Lambda^{\bullet} V(n))
       = \YM(n) \otimes \Lambda^{\bullet \geq p} V(n).     
\end{equation}
It may be easily verified that $F_{p} (\YM(n) \otimes \Lambda^{\bullet} V(n))$ is decreasing, bounded, multiplicative and compatible with differentials,
i.e. $d(F_{p} (\YM(n) \otimes \Lambda^{\bullet} V(n))) \subseteq F_{p} (\YM(n) \otimes \Lambda^{\bullet} V(n))$.
Furthermore,
\[     d(F_{p} (\YM(n) \otimes \Lambda^{\bullet} V(n)))
       \subseteq F_{p+1} (\YM(n) \otimes \Lambda^{\bullet} V(n)).     \]
Notice that $\epsilon (F_{p}(\YM(n) \otimes \Lambda^{\bullet} V(n))) = 0$, if $p \geq 1$.

As a consequence, the associated graded algebra to this filtration $\mathrm{Gr}_{F_{\bullet}}(\YM(n) \otimes \Lambda^{\bullet} V(n))$
is an augmented differential graded algebra, provided with zero differential and augmentation induced by $\epsilon$.

As before, we will consider $\YM(n)$ as an augmented differential graded algebra concentrated in degree zero and augmentation
$\epsilon_{\ym(n)}$, and $\Lambda^{\bullet} V(n)$ as an augmented differential graded algebra with the usual grading (i.e. given by $\bullet$)
and the usual augmentation, both with zero differential.
The associated graded algebra $\mathrm{Gr}_{F_{\bullet}}(\YM(n) \otimes \Lambda^{\bullet} V(n))$ 
is the tensor product (in the category of augmented differential graded algebras) of the algebras $\YM(n)$ and $\Lambda^{\bullet} V(n)$.

The filtration \eqref{eq:anotherfiltration} induces a decreasing and bounded above filtration of coaugmented coalgebras on
\[     B^{+}(\YM(n) \otimes \Lambda^{\bullet} V(n)),     \]
which we denote $F_{\bullet} (B^{+}(\YM(n) \otimes \Lambda^{\bullet} V(n)))$ (cf. \cite{Lef1}, Sec. 1.3.2).
By its very definition,
\[     \mathrm{Gr}_{F_{\bullet}}(B^{+}(\YM(n) \otimes \Lambda^{\bullet} V(n)))
       = B^{+} (\mathrm{Gr}_{F_{\bullet}}(\YM(n) \otimes \Lambda^{\bullet} V(n))).     \]

Since the associated graded algebra $\mathrm{Gr}_{F_{\bullet}}(\YM(n) \otimes \Lambda^{\bullet} V(n))$ is the tensor product
of the algebras $\YM(n)$ and $\Lambda^{\bullet} V(n)$, then 
\[     B^{+} (\mathrm{Gr}_{F_{\bullet}}(\YM(n) \otimes \Lambda^{\bullet} V(n)))
       \simeq_{q} B^{+} (\YM(n)) \otimes B^{+} (\Lambda^{\bullet} V(n)).     \]
Again using that $B^{+}(\YM(n))$ is quasi-isomorphic as a graded $k$-vector space to the tensor product over $\YM(n)^{e}$
of any projective resolution of graded $\YM(n)^{e}$-modules of $k$ with $k$, we obtain that $B^{+}(\YM(n)) \simeq_{q} C_{\bullet} (\YM(n),k)$.

Besides, since $\Lambda^{\bullet} V(n)$ is Koszul, there is a quasi-isomorphism $B^{+}(\Lambda^{\bullet} V(n)) \simeq_{q} S^{\bullet}(V(n))$,
where we consider $S^{\bullet}(V(n))$ as a differential graded algebra provided with the usual grading (i.e. given by $\bullet$)
and zero differential (cf. \cite{PP1}, Sec. 2.1, Example; Sec. 2.3).

Finally,
\[     \mathrm{Gr}_{F_{\bullet}}(B^{+}(\YM(n) \otimes \Lambda^{\bullet} V(n)))
       \simeq_{q} C_{\bullet} (\YM(n),k) \otimes S(V(n)).     \]

Notice that the space on the right is the zero step of the weakly convergent spectral sequence in Corollary \ref{coro:sucesionespectral}.

Moreover, the spectral sequence associated to the filtration $F_{\bullet}(B^{+}(\YM(n) \otimes \Lambda^{\bullet} V(n)))$ is convergent, for
the filtration being decreasing, exhaustive, Hausdorff and bounded above.
From the Comparison Theorem for spectral sequences (cf. \cite{Wei1}, Thm. 5.2.12) and the Mapping Lemma for 
$E^{\infty}$ (cf. \cite{Wei1}, Exercise 5.2.3), it follows that
\[     H_{\bullet} (B^{+}(\YM(n) \otimes \Lambda^{\bullet} V(n))) \simeq H_{\bullet} (C_{\bullet} (\YM(n),S(V(n)))).     \]
By Propositi\'on \ref{prop:quasiiso}
\begin{equation}
\label{eq:isogen}
     H_{\bullet} (\tym(n),k) \simeq H_{\bullet} (\TYM(n),k)) \simeq H_{\bullet} (B^{+}(\TYM(n)))
                               \simeq H_{\bullet} (C_{\bullet} (\YM(n),S(V(n)))).
\end{equation}

We now recall some useful results concerning cohomology of algebras. 
\begin{prop}
\label{prop:algebraLielibre}
Let $k$ be a commutative ring with unit and let $\g$ be a Lie algebra over $k$. 
If $H^{\bullet}(\g, M) = 0$, for all $\g$-module $M$ and $\bullet \geq 2$.
Then $\g$ is a free Lie algebra.

Let $W$ be a $k$-module.
If $M$ is a bimodule over the free algebra $T(W)$, then $H_{\bullet}(T(W), M) = 0$ and $H^{\bullet}(T(W), M) = 0$, for all $\bullet \geq 2$.
\end{prop}
\noindent\textbf{Proof.}
Cf. \cite{Wei1}, Prop. 9.1.6. and Ex. 7.6.3.
\qed

\begin{rem}
\label{obs:librelielibrealgebra}
From Proposition \ref{prop:algebraLielibre}, we immediately see that if the universal enveloping algebra $\U(\g)$ of
a Lie algebra $\g$ is free, then $\g$ is also free.
The converse is also clear.
\end{rem}

\begin{teo}
\label{teo:bergerlibregen}
Let $W$ be a graded $k$-vector space concentrated in degree $1$,
let $I$ be an homogeneous ideal generated in degrees greater or equal to $2$ and let $A = T(W)/I$.
We may write $I = \bigoplus_{m \in \NN_{\geq 2}} I_{m}$, where $I_{m}$ denotes the $m$-th homogeneous component.
Then, for each $m \geq 2$, it is possible to choose a $k$-linear subspace $R_{m} \subseteq I_{m}$ such that
\begin{align*}
   I_{2} &= R_{2},
   \\
   I_{m} &= R_{m} \oplus \left( \sum_{i+j+l=m, 2 \leq j < m} W^{\otimes i} \otimes R_{j} \otimes W^{\otimes l} \right).
\end{align*}
Taking $R = \bigoplus_{m \in \NN_{\geq 2}} R_{m}$, the homology is
\begin{align*}
   H_{0} (A,k) &= \Tor_{0}^{A} (k,k) = k,
   \\
   H_{1} (A,k) &= \Tor_{1}^{A} (k,k) = V,
   \\
   H_{2} (A,k) &= \Tor_{2}^{A} (k,k) = R.
\end{align*}
\end{teo}
\noindent\textbf{Proof.}
Cf. \cite{Ber2}, p. 9 and \cite{HKL1}, Proof of the Prop. in the Appendix. 
\qed

Using Proposition \ref{prop:homologiayangmillssimetrica}, we obtain that $H_{2} (\TYM(n) , k) = 0$.
If we consider the vector space $V(n)$ as a graded vector space concentrated in degree $2$, the algebra $\TYM(n) = \U(\tym(n))$
is graded, and Theorem \ref{teo:bergerlibregen} implies that $\TYM(n)$ is a free graded algebra,
so a fortiori, $\TYM(n)$ is a free algebra, when we forget the grading.
Furthermore, from Remark \ref{obs:librelielibrealgebra} it yields that $\tym(n)$ is a free Lie algebra, in the graded sense or not, due to
Remark \ref{obs:libregraduado}.

We have proved the following theorem, which is a key result in order to study representations of the Yang-Mills algebra.
\begin{teo}
\label{teo:tyangmillslibre}
The Lie algebra $\tym(n)$ is free and the same holds for the associative algebra $\TYM(n)$.
\end{teo}

The last paragraphs of this section are devoted to a description of the graded vector space which generates the above free algebra.

Let us call $W(n)$ the graded $k$-vector space satisfying $\f_{gr}(W(n)) \simeq \tym(n)$, or equivalently $T_{gr}(W(n)) \simeq \TYM(n)$, where we are
using the special grading.
By Theorem \ref{teo:bergerlibregen}, we see that $H_{1}(\TYM(n),k) \simeq W(n)$ as graded vector spaces, since all morphisms considered in this section are homogeneous of degree $0$.

Using that, by the isomorphism \eqref{eq:isogen}, $H_{\bullet} (\TYM(n),k)) \simeq H_{\bullet} (\ym(n),S(V(n)))$ and that
Proposition \ref{prop:homologiayangmillssimetrica} implies the following isomorphism of graded vector spaces
\[     H_{1} (\ym(n),S(V(n))) \simeq \sum_{p \in \NN} H_{1}^{p} (\ym(n),S(V(n)))     \]
where $H_{1}^{p} (\ym(n),S(V(n)))$ lives in degree $p+1$,
we immediately see that the homogeneous component $W(n)_{m}$ of degree $m$ of $W(n)$ is isomorphic to $H_{1}^{m-1} (\ym(n),S(V(n)))$.

In particular, $W(2) \simeq H_{1}(\ym(2),S(V(2))) = H_{1}^{1} (\ym(2),S(V(2))) \simeq k[-2]$ (as $k$-vector spaces) and $W(n)$ is infinite dimensional if $n \geq 3$.
For instance, it follows easily from Proposition \ref{prop:homologiayangmillssimetrica} that for $n=3$,
\begin{equation}
\label{eq:W(3)}
   W(3) \simeq H_{1}(\ym(3),S(V(3))) = \bigoplus_{p \in \NN_{0}} H_{1}^{p} (\ym(3),S(V(3))),
\end{equation}
where $\dim_{k}(H_{1}^{p} (\ym(3),S(V(3)))) = 2 p + 1$ for $p \in \NN$ and $H_{1}^{0} (\ym(3),S(V(3)))) = 0$.

\medskip

We recall the following definition.
\begin{defi}(cf. \cite{PP1}, \cite{La2}, Ch. X, \S 6)
The \emph{Hilbert series} of a $\ZZ$-graded $k$-vector space $V = \bigoplus_{n \in \ZZ} V_{n}$ (also called \emph{Poincar\'e series}) is
the formal power series in $\ZZ[[t, t^{-1}]]$ given by
\[     h_{V} (t) = \sum_{n \in \ZZ} \dim_{k}(V_{n}) t^{n}.     \]
We will also denote it simply $V(t)$.
\end{defi}

We will now compute the Hilbert series $W(n)(t)$ of $W(n)$ provided with the usual grading.
In order to do so, we must notice that if
\[     0 \rightarrow M' \overset{f}{\rightarrow} M \overset{g}{\rightarrow} M'' \rightarrow 0     \]
is a short exact sequence of graded vector spaces with homogeneous morphisms $f$ and $g$ of degree $0$, then
\[     M(t) = M''(t) + M'(t).     \]

The previous identity implies that that the Hilbert series is an \emph{Euler-Poincar\'e map} (cf. \cite{La2}, Ch. III, \S 8)
when considering the category of graded vector spaces with homogeneous morphisms of degree $0$.
As a consequence, if we consider the \emph{Euler-Poincar\'e characteristic} of a complex $(C_{\bullet},d_{\bullet})$
of graded vector spaces with homogeneous morphisms, which is defined as
\[     \chi(C)(t) = \sum_{i \in \ZZ} (-1)^{i} C_{i}(t),     \]
it turns out that it coincides with the Euler-Poincar\'e characteristic of its homology
(cf. \cite{La2}, Ch. XX, \S 3, Thm. 3.1).

On the one hand, the Euler-Poincar\'e characteristic of the complex \eqref{eq:complejocohomologiayangmillssimetrica} is
\[     \chi_{C_{\bullet}(\YM(n),S(V(n)))(t)} = S(V(n))(t) - n t S(V(n))(t) + n t^{3} S(V(n))(t) - t^{4} S(V(n))(t)
                                                          =  \frac{1- n t + n t^{3} - t^{4}}{(1-t)^{n}},     \]
since $S(V(n))(t) = (1-t)^{-n}$.

On the other hand, since the morphisms of the complex \eqref{eq:complejocohomologiayangmillssimetrica}
are homogeneous of degree $0$, its Euler-Poincar\'e characteristic coincides with the one of its homology.

Since $H_{0}(\ym(n), S(V(n))) \simeq k$, $H_{1}(\ym(n), S(V(n))) \simeq W(n)$ and $H_{2}(\ym(n), S(V(n))) = H_{3}(\ym(n), S(V(n))) = 0$ as graded vector spaces, the Euler-Poincar\'e characteristic of the homology is
\[     \chi_{H_{\bullet}(C_{\bullet}(\YM(n),S(V(n))))(t)} = 1 - W(n)(t).     \]
Finally,
\[     W(n)(t) = \frac{(1-t)^{n}-1+n t -n t^{3} +t^{4}}{(1-t)^{n}}.     \]

We have thus proved the following proposition.
\begin{prop}
\label{prop:serieHilbertW(n)}
The Hilbert series of the graded vector space $W(n)$ with the usual grading is given by
\[     W(n)(t) = \frac{(1-t)^{n}-1+n t -n t^{3} +t^{4}}{(1-t)^{n}}.     \]
\end{prop}

It is trivially verified that the term of degree $m$ of $W(n)(t)$ coincides with the computation of $\dim_{k}(W(n)_{m})$
obtained from Proposition \ref{prop:homologiayangmillssimetrica}.

Observe that for $n = 2$, the Yang-Mills algebra $\YM(2)$ is noetherian.
This is a direct consequence of the isomorphism $\YM(2) \simeq \U(\h_{1})$ (cf. Example \ref{ejem:heisenberg}).
Since $\h_{1}$ is finite dimensional, then $S(\h_{1})$ is noetherian, so $\U(\h_{1})$ is also noetherian, because
its associated graded algebra $S(\h_{1})$ is (cf. \cite{Dix1}, Corollary 2.3.8, p. 76).

However, for $n \geq 3$, Yang-Mills algebras $\YM(n)$ are non noetherian.
In order to prove this fact we proceed as follows.

On the one side, $\YM(n) \supset \TYM(n)$ is a (left and right) free extension of algebras,
which is not finite but finitely generated (cf. \cite{Wei1}, Coro. 7.3.9).
One set of generators of the extension $\YM(n) \supset \TYM(n)$ is $\{ x_{1}, \dots, x_{n} \}$.

On the other side, since $\TYM(n)$ is a free algebra with an infinite set of generators if $n \geq 3$ (cf. Thm. \ref{teo:tyangmillslibre}, 
Prop. \ref{prop:serieHilbertW(n)}), it is not noetherian, which implies that $\YM(n)$ is not noetherian.
This is a direct consequence of the following simple fact:
If $A \supset B$ is an extension of algebras such that $A$ is a right (resp. left) free $B$-module and $I \triangleleft B$ is a left (resp. right) ideal,
then $A.I$ is a left (resp. right) ideal of $A$ that satisfies that $A.I \cap B = I$.
This property directly yields that if $B$ is not left (resp. right) noetherian, then $A$ is not left (resp. right) noetherian.

Let us thus show the previously stated property for free extensions of algebras.
We shall prove it in the case that $A$ is a right free $B$-module and $I \triangleleft B$ is a left ideal, the other being analogous.

The inclusion $I \subseteq A.I \cap B$ is clear.
We will prove the other one.
Let $\B = \{ a_{j} \}_{j \in J}$ be a basis of $A$ as a right $B$-module.
We assume without loss of generality that $1 = a_{j_{0}} \in \B$.
Since every element of $A$ can be written as $\sum_{j \in J} a_{j} b_{j}$, with $b_{j} \in B$,
an element $x \in A.I \cap B$ may be written as $x = \sum_{j \in J} a_{J} c_{j}$, where $c_{j} \in I$.
But $x \in B$, so $c_{j_{0}} = x$ and $a_{j} = 0$, for all $j \in J$, $j \neq j_{0}$.
Hence, $x = a_{j_{0}} c_{j_{0}} = 1 c_{j_{0}} = c_{j_{0}} \in I$.

\begin{rem}
\label{rem:nonfinitedimensional}
The fact that $\YM(n) = \U(\ym(n))$ is not noetherian for $n \geq 3$ implies that $\ym(n)$ is not finite dimensional for $n \geq 3$.
\end{rem}

\section{Main theorem: relation between Yang-Mills algebras and Weyl algebras}
\label{sec:main}

The aim of this last section is to prove that all the Weyl algebras $A_{r}(k)$ ($r \in \NN$) are epimorphic images of all Yang-Mills algebras
$\YM(n)$
for $n \geq 3$.
In order to do so we make intensive use of $\TYM(n)$.
As a consequence, the representations of all $A_{r}(k)$ are also representations of $\YM(n)$ ($r \in \NN$, $n \geq 3$).
These families of representations have been previously studied by Bavula and Bekkert in \cite{BB1} and are enough to separate points of
$\YM(n)$.

Since $\YM(3)$ is a quotient of $\YM(n)$, as an algebra, for every $n \geq 3$, it will be sufficient to prove that that the Weyl algebras $A_{r}(k)$
are epimorphic images of $\YM(3)$.
Our first step is to give explicit bases for the quotients $\ym(3)/\C^{j}(\ym(3))$, for $j= 1,2,3,4$.
The elements of these bases will be useful while defining the epimorphisms onto $A_{r}(k)$.

Using the fact that $\YM(n)$ is Koszul, Connes and Dubois-Violette obtain the Hilbert series of the Yang-Mills algebra $\YM(n)$
(cf. \cite{CD1}, Corollary 3)
\[     h_{\YM(n)}(t) = \frac{1}{(1-t^{2})(1-n t + t^{2})}.     \]

Using the Poincar\'e-Birkhoff-Witt Theorem it is possible to deduce the dimensions of the homogeneous spaces of the Yang-Mills algebra
$N(n)_{j} = \dim_{k}(\ym(n)_{j})$ ($j \in \NN$) from the equality
\[     \prod_{j \in \NN} \left(\frac{1}{1-t^{j}}\right)^{N(n)_{j}} = h_{\YM(n)}(t).     \]
Connes and Dubois-Violette find in \cite{CD1} the direct formula ($j \geq 3$)
\[     N(n)_{j} = \frac{1}{j} \sum_{k=1}^{j} \mu(\frac{j}{k}) (t_{1}^{k}+t_{2}^{k}),     \]
where $t_{1}$ y $t_{2}$ are the roots of the polynomial $t^{2} - n t +1 = 0$ and $\mu(x)$ is the M\"obius function.

For $\ym(3)$, the sequence of dimensions $N(3)_{j}$ ($j \in \NN$) is (cf. \cite{Sl1}, sequence A$072337$)
\[     3, 3, 5, 10, 24, 50, 120, 270, 640, 1500, 3600, 8610, 20880, 50700, 124024, 304290, 750120, \dots     \]

An ordered basis for the quotient algebra $\ym(3)/\C^{1}(\ym(3))$ of the Yang-Mills algebra is given by
\[     \B_{1} = \{ x_{1}, x_{2}, x_{3} \}.      \]

For the quotient $\ym(3)/\C^{2}(\ym(3))$, a possible ordered basis is
\[     \B_{2} = \{ x_{1}, x_{2}, x_{3}, x_{12}, x_{13} , x_{23} \},     \]
where $x_{ij} = [x_{i} , x_{j}]$, ($i, j = 1, 2, 3$).
To prove that it is indeed a basis we must only show that it generates $\ym(3)/\C^{2}(\ym(3))$ (since $\#(\B_{2}) = 6$),
which is obtained using that
$x_{ij} = - x_{ji}$.

The following set is a basis of $\ym(3)/\C^{3}(\ym(3))$:
\[     \B_{3} = \{ x_{1}, x_{2}, x_{3}, x_{12}, x_{13} , x_{23}, x_{112}, x_{221}, x_{113}, x_{123}, x_{312} \},     \]
where we denote $x_{ijk} = [x_{i},[x_{j},x_{k}]]$.
We also define $J_{3} = \{ (112), (221), (113), (123), (312) \}$ the set of triple indices of $\B_{3}$.

Let us prove that $\B_{3}$ is a basis.
As before, we only have to prove that it generates $\ym(3)/\C^{3}(\ym(3))$.
This is direct, as we can see from the Yang-Mills relations
\[     x_{332} = - x_{112}, \hskip 0.5cm x_{331} = - x_{221}, \hskip 0.5cm x_{223} = - x_{113},    \]
and relations given by antisymmetry and Jacobi identity, i.e. $x_{ijk} = - x_{ikj}$, and $x_{213} = x_{123} + x_{312}$.

The case $\ym(3)/\C^{4}(\ym(3))$ is a little more complicated.
We shall prove that
\begin{align*}
   \B_{4} = \{ &x_{1}, x_{2}, x_{3}, x_{12}, x_{13} , x_{23}, x_{112}, x_{221}, x_{113}, x_{123}, x_{312},
   \\
   &x_{1112}, x_{1221}, x_{1113}, x_{1123}, x_{2221}, x_{2113}, x_{2312}, x_{3112}, x_{3221}, x_{3312} \},
\end{align*}
is an ordered basis, where $x_{ijkl} = [x_{i},[x_{j},[x_{k},x_{l}]]]$.
Analogously, we define the set of indices of $\B_{4}$
\[     J_{4} = \{ (1112), (1221), (1113), (1123), (2221), (2113), (2312), (3112), (3221), (3312) \}.     \]

In order to prove that $\B_{4}$ is a basis it suffices again to verify that it generates $\ym(3)/\C^{4}(\ym(3))$.
On the one hand, taking into account that
\[     [[x_{i},x_{j}],[x_{k},x_{l}]] = [[[x_{i},x_{j}],x_{k}],x_{l}] + [x_{k},[[x_{i},x_{j}],x_{l}]]
                                     = [x_{l},[x_{k},[x_{i},x_{j}]]] + [x_{k},[x_{l},[x_{j},x_{i}]]]
                                     = x_{lkij} + x_{klji}
\]
and $[[[x_{i},x_{j}],x_{k}],x_{l}] = [x_{l},[x_{k},[x_{i},x_{j}]]] = x_{lkij}$, 
the set $\{ x_{ijkl} : i, j, k, l = 1, 2, 3 \} \cup \B_{3}$ is a system of generators.
We shall prove that it is generated by $\B_{4}$.
In fact, we only need to prove that it generates the set
\[     \{  x_{i112}, x_{i221}, x_{i113}, x_{i123}, x_{i312} : i = 1, 2, 3 \},     \]
because $\{ x_{112}, x_{221}, x_{113}, x_{123}, x_{312} \}$ are generators of the homogeneous elements of degree $3$.
This last statement is direct:
\begin{align*}
     x_{3113} &= x_{1221}, \hskip 0.5cm x_{2112} = - x_{1221}, \hskip 0.5cm
     x_{2123} = x_{3221} + x_{2312} - x_{1113},
     \\
     x_{1312} &= \frac{x_{3112} + x_{2113} - x_{1123}}{2}, \hskip 0.95cm x_{3123} = \frac{x_{1112}+x_{2221}-x_{3312}}{2}.
\end{align*}

We shall now briefly recall a version of the Kirillov orbit method by J. Dixmier, which we will employ. 
We first recall that a bilateral ideal $I \triangleleft A$ of an algebra $A$ is called \emph{prime} if $I \neq A$ 
and if $J, K \triangleleft A/I$ are two nonzero bilateral ideals of the quotient algebra $A/I$, then $J K \neq \{0\}$. 
We say that $I \triangleleft A$ is \emph{completely prime} if $A/I$ is a domain.
Observe that every completely prime ideal is prime (cf. \cite{Dix1}, 3.1.6).
A bilateral ideal $I \triangleleft A$ is called \emph{semiprime} if $I \neq A$ and every nilpotent bilateral ideal $J \triangleleft A/I$ is zero.
Note that an intersection of semiprime ideals is semiprime and every prime ideal is semiprime (cf. \cite{Dix1}, 3.1.6).

On the other hand, a bilateral ideal $I \triangleleft A$ is called \emph{primitive} if it the annihilating ideal of a simple left $A$-module, 
and it called \emph{maximal} if $I \neq A$ and if it is maximal in the lattice of bilateral ideals of $A$, ordered by inclusion.
Every maximal ideal is primitive (cf. \cite{Dix1}, 3.1.6).

If $\g$ is finite dimensional Lie algebra, $\U(\g)$ is a noetherian domain, so it has a skew-field of fractions $\Frac(\U(\g))$ 
(cf. \cite{Dix1}, 3.1.16 and Thm. 3.6.12).
Let $I \triangleleft \U(\g)$ be a semiprime ideal, so $\U(\g)/I$ also has a skew-field of fractions $\Frac(\U(\g)/I)$.
We say that $I$ is \emph{rational} if $\Z(\Frac(\U(\g)/I)) = k$.
Every rational ideal of $\U(\g)$ is primitive and, $k$ being algebraically closed, the converse also holds (cf. \cite{Dix1}, Thm. 4.5.7).

We have the following proposition. 
\begin{prop}
\label{prop:lieweyl}
Let $I \triangleleft \U(\g)$ be a bilateral ideal of the universal enveloping algebra of a nilpotent Lie algebra $\g$ of finite dimension.
The following are equivalent:
\begin{itemize}
\item[(i)] $I$ is primitive.
\item[(ii)] $I$ is maximal.
\item[(iii)] There exists $r \in \NN$ such that $\U(\g)/I \simeq A_{r}(k)$.
\item[(iv)] $I$ is the kernel of a simple representation of $\U(\g)$.
\end{itemize}
\end{prop}
\noindent\textbf{Proof.} 
Cf. \cite{Dix1}, Prop. 4.7.4, Thm. 4.7.9.
\qed

If $I \triangleleft \U(\g)$ is a bilateral ideal satisfying either of the previous equivalent conditions, the positive integer $r$ 
(uniquely determined) such that $\U(\g)/I \simeq A_{r}(k)$ is called the \emph{weight} of the ideal $I$ (cf. \cite{Dix1}, 4.7.10).

Let us suppose that $\g$ is nilpotent Lie algebra of finite dimension. 
Given $f \in \g^{*}$, a \emph{polarization} of $f$ is a subalgebra $\h \subset \g$
such that it is \emph{subordinated} to $f$, i.e. $f([\h,\h]) = 0$, and it is maximal with respect to the previous property
(cf. \cite{Dix1}, Section 1.12). 
There exists a canonical way to construct polarizations, which are called \emph{standard}, of a linear form $f \in \g^{*}$ 
(cf. \cite{Dix1}, Prop. 1.12.10). 
It is easily verified that the dimension of a polarization $\h$ of $f$ must be (cf. \cite{Dix1}, 1.12.1)
\begin{equation}
\label{eq:dim}
     \frac{\dim_{k}(\g) + \dim_{k}(\g^{f})}{2}.
\end{equation}
Furthermore, $\h$ is a polarization of $f$ on $\g$ if and only if $\h$ is a subalgebra subordinated to $f$ of dimension
$(\dim_{k}(\g) + \dim_{k}(\g^{f}))/2$ (cf. \cite{Dix1}, 1.12.8). 
The weight of a primitive ideal $I(f)$ is given by $r = \dim_{k}(\g/\g^{f})/2$ (cf. \cite{Dix1}, Prop. 6.2.2), or using
identity \eqref{eq:dim}, $r = \dim_{k}(\g/\h_{f})$, where $\h_{f}$ is a polarization of $f$.

We shall now explain the connection between rational ideals and polarizations. 
If $f \in \g^{*}$ be a linear functional and $\h_{f}$ a polarization of $f$,
we may define a representation of $\h_{f}$ on the vector space $k.v_{f}$ of dimension $1$ by means of 
$x.v_{f} = (f(x) + \tr_{\g/\h_{f}} (\ad_{\g}x))v_{f}$, for $x \in \h_{f}$ and $\tr_{\g/\h_{f}} = \tr_{\g} - \tr_{\h_{f}}$.
Therefore, we can consider the induced $\U(\g)$-module $V_{f} = \U(\g) \otimes_{\U(\h_{f})} k.v_{f}$.
If we denote the corresponding action $\rho : \U(\g) \rightarrow \End_{k}(V_{f})$, $I(f) = \Ker(\rho)$ is a bilateral ideal of the 
enveloping algebra $\U(\g)$.

In the previous notation we have omitted the polarization in $I(f)$. 
This is justified by the following proposition, which states even more. 
\begin{prop}
Let $\g$ be a nilpotent Lie algebra of finite dimension, let $f \in \g^{*}$ and let $\h_{f}$ and $\h'_{f}$ be two polarizations of $f$.
If we denote $\rho : \U(\g) \rightarrow \End_{k}(V_{f})$ and $\rho' : \U(\g) \rightarrow \End_{k}(V'_{f})$ 
the corresponding representations constructed following the previous method, then $\Ker(\rho) = \Ker(\rho')$. 
This ideal is primitive. 
 
On the other hand, if $I$ is a primitive ideal of $\U(\g)$, then there exists $f \in \g^{*}$ such that $I = I(f)$.
\end{prop}
\noindent\textbf{Proof.} 
Cf. \cite{Dix1}, Thm. 6.1.1, Thm. 6.1.4 and Thm. 6.1.7. 
\qed

The group $\Aut(\g)$ is an algebraic group whose associated Lie algebra is $\Der(\g)$.
Let $\mathfrak{a}$ denote the algebraic Lie algebra generated by the ideal $\InnDer(\g)$ in $\Der(\g)$.
The irreducible algebraic group $G$ associated to $\mathfrak{a}$ is called the \emph{adjoint algebraic group} of $\g$.
It is a subgroup of $\Aut(\g)$.
If $\mathfrak{a} = \InnDer(\g)$, $G$ is called the \emph{adjoint group} of $\g$.
As a consequence, the group $G$ acts on the Lie algebra $\g$, so it also acts on $\g^{*}$ with the dual action, which is called \emph{coadjoint}. 
\begin{teo}
Let $\g$ be a nilpotent Lie algebra of finite dimension and let $f$ and $f'$ be two linear forms on $\g$. 
If $I(f)$ and $I(f')$ are the corresponding bilateral ideals of $\U(\g)$, then $I(f) = I(f')$ if and only if there is $g \in G$ such that 
$f = g.f'$.
\end{teo}
\noindent\textbf{Proof.} 
Cf. \cite{Dix1}, Prop. 6.2.3.
\qed

The previous results imply that, for a nilpotent Lie algebra of finite dimension there exists a bijection
\[     I : \g^{*}/G \rightarrow \Prim (\U(\g))     \]
between the set of classes of linear forms on $\g$ under the coadjoint action and the set of primitive ideals of $\U(\g)$.

We shall now use the previous method for the first quotients of the Yang-Mills algebra $\ym(n)$ by the the ideals of the low derived sequence 
$\C^{l}(\ym(n))$ in order to find the possible weights for each step.  

Given $f \in (\ym(3)/\C^{2}(\ym(3)))^{*}$,
\[     f = \sum_{i=1}^{3} c_{i} x_{i}^{*} + \sum_{i j} c_{ij} x_{ij}^{*},     \]
with $c_{12} = 0$, $c_{13} \neq 0$ and $c_{23} \neq 0$, we obtain a polarization associated to $f$ as follows.
Consider
\[     \h_{f} = k.x_{1} \oplus k.x_{2} \oplus k.x_{12} \oplus k.x_{13} \oplus k.x_{23}.     \]
Since the weight of $I(f)$ (i.e. the positive integer $n$ such that $\U(\g)/I(f) \simeq A_{n}(k)$)
is equal to $\dim_{k} (\g/\h_{f})$ (cf. \cite{Dix1}, 1.2.1, 1.2.8 and Prop. 6.2.2),
we obtain that $\U(\ym(3)/\C^{2}(\ym(3)))/I(f) \simeq A_{1}(k)$.
It is easy to show that there are no higher weights for $\ym(3)/\C^{2}(\ym(3))$.

Any element of $(\ym(3)/\C^{3}(\ym(3)))^{*}$ may be written as
\[     f = \sum_{i=1}^{3} c_{i} x_{i}^{*} + \sum_{i < j} c_{ij} x_{ij}^{*} + \sum_{(ijk) \in J_{3}} c_{ijk} x_{ijk}^{*}.     \]
If $c_{112} = c_{123} = 1$, $c_{113} = c_{221} = c_{312} = 0$, we find a polarization associated to $f$ as follows:
\[     \h_{f} = k.x_{12} \oplus k.x_{13} \oplus k.x_{23} \oplus k.x_{112} \oplus k.x_{221} \oplus k.x_{113}
                                \oplus k.x_{123} \oplus k.x_{312},     \]
and hence $\U(\ym(3)/\C^{3}(\ym(3)))/I(f) \simeq A_{3}(k)$.
If $c_{112} = 1$, $c_{221} = c_{113} = c_{312} = c_{123} = 0$,
\[     \h_{f} = k.x_{2} \oplus k.x_{12} \oplus k.x_{13} \oplus k.x_{23} \oplus k.x_{112} \oplus k.x_{221} \oplus k.x_{113}
                                \oplus k.x_{123} \oplus k.x_{312},     \]
and hence $\U(\ym(3)/\C^{3}(\ym(3)))/I(f) \simeq A_{2}(k)$.

For $\ym(3)/\C^{4}(\ym(3))$, every linear functional has the form
\[     f = \sum_{i=1}^{3} c_{i} x_{i}^{*} + \sum_{i < j} c_{ij} x_{ij}^{*} + \sum_{(ijk) \in J_{3}} c_{ijk} x_{ijk}^{*}
       + \sum_{(ijkl) \in J_{4}} c_{ijkl} x_{ijkl}^{*}.     \]
If $c_{312} = c_{2312} = c_{1112} = 1$, and all other coefficients are zero, so, taking
\[     \h_{f} = k.x_{12} \oplus k.x_{13} \oplus k.x_{112} \oplus k.x_{221} \oplus k.x_{113} \oplus k.x_{123}
                   \oplus k.x_{312} \oplus \bigoplus_{(ijkl) \in J_{4}} k.x_{ijkl},     \]
we obtain that $\U(\ym(3)/\C^{4}(\ym(3)))/I(f) \simeq A_{4}(k)$.

As a consequence, we have just proved that, given $n$ such that $1 \leq n \leq 4$, there exists an ideal $I$ in $\U(\ym(3))$ satisfying
$A_{n}(k) \simeq \U(\ym(3))/I$.
In fact, a stronger statement is true:

\begin{teo}
\label{teo:yangmillsweyl}
Let $r \in \NN$ be a positive integer. There exists a surjective homomorphism of algebras
\[     \U(\ym(3)) \twoheadrightarrow A_{r}(k).     \]
Furthermore, there exists $l \in \NN$ such that we can choose this homomorphism satisfying that it factors through the quotient
$\U(\ym(3)/\mathcal{C}^{l}(\ym(3)))$
\[
\xymatrix@C-20pt@R-10pt { \U(\ym(3)) \ar@{->>}[rr] \ar@{->>}[rd] &
& A_{r}(k)
\\
& \U(\ym(3)/\mathcal{C}^{l}(\ym(3))) \ar@{->>}[ru] & }
\]
\end{teo}

Before giving the proof of this theorem, we shall state the next corollary, our main result, which follows readily from
Theorem \ref{teo:yangmillsweyl} and the fact that that every $\U(\ym(3))$ is a quotient of every $\U(\ym(n))$ for $n \geq 4$.
\begin{coro}
\label{coro:yangmillsweyl}
Let $r, n \in \NN$ be two positive integers, satisfying $n \geq 3$.
There exists a surjective homomorphism of $k$-algebras
\[       \U(\ym(n)) \twoheadrightarrow A_{r}(k).     \]
Furthermore, there exists $l \in \NN$ such that we can choose this morphism in such a way that it factors through the quotient
$\U(\ym(n)/\mathcal{C}^{l}(\ym(n)))$
\[
\xymatrix@C-20pt@R-10pt { \U(\ym(n)) \ar@{->>}[rr] \ar@{->>}[rd] &
& A_{r}(k)
\\
& \U(\ym(n)/\mathcal{C}^{l}(\ym(n))) \ar@{->>}[ru] & }
\]
\end{coro}

The previous corollary shows the strong link between representations of the Yang-Mills algebra $\YM(n)$
and representations of the Weyl algebra $A_{r}(k)$.
Concretely, using the characterization of the category of representations over the latter exhibited in \cite{BB1},
we are able to give a description of a subcategory of infinite dimensional modules over the Yang-Mills algebra $\YM(n)$, for $n \geq 3$,
by means of the surjection in Corollary \ref{coro:yangmillsweyl}.
In order to do so, we shall briefly recall the results presented in \cite{BB1} and we refer to it for the general definitions of
weight modules, generalized weight modules and orbits.
Given an orbit $\O$, we denote by $\W(\O)$ and $\GW(\O)$ the subcategories of the categories of modules over $A_{r}(k)$
consisting of weight modules and generalized modules with support in the orbit $\O$, respectively.

Also, given an orbit $\O$, Bavula and Bekkert construct certain categories $\L_{\O,\mathrm{nil}}$ and $\L_{\O,1}$, whose set of objects is $\O$,
isomorphic to $\GW(\O)$ and $\W(\O)$, respectively.
Furthermore, if $\O$ is an orbit and $M$ is a $\L_{\O, \mathrm{nil}}$-module (resp. $\L_{\O,1}$-module), they define an $A_{r}(k)$-module $M_{\O}$.

Depending on the orbit, either $\L_{\O,\mathrm{nil}}$ or $\L_{\O,1}$ may be isomorphic to one of the tame categories exhibited in Figure
\ref{fig:tame} (cf. \cite{BB1}, Sec. 2.5 to 2.7):

\begin{center}
\begin{figure}[H]
\[
\xymatrix@R-10pt
{
\A = k,
&
&
&
&
&
\\
\B = k\cl{t},
&
&
&
&
&
\text{$t$ nilpotent,}
\\
\C = k\cl{t, t^{-1}},
&
&
&
&
&
\\
\D
&
\bullet^{1}
\ar@<1ex>[r]^{a}
&
\bullet^{2}
\ar@<1ex>[l]^{b}
&
&
&
\text{$a b$ nilpotent,}
\\
\E
&
\bullet^{1}
\ar@(ul,dl)_{a_{11}}
\ar@<1ex>[r]^{a_{21}}
&
\bullet^{2}
\ar@(ur,dr)^{a_{22}}
\ar@<1ex>[l]^{a_{12}}
&
&
&
\text{$a_{11}^{2} = a_{12} a_{21}$, $a_{22}^{2} = a_{21}a_{12}$, $a_{21}a_{11} = a_{22}a_{21}$,}
\\
&
&
&
&
&
\text{$a_{12}a_{22} = a_{11}a_{12}$; $a_{11}$ and $a_{22}$ nilpotent,}
\\
\F
&
\bullet^{1}
\ar@<1ex>[r]^{a_{2}}
\ar@<1ex>[d]^{b_{1}}
&
\bullet^{2}
\ar@<1ex>[l]^{b_{2}}
\ar@<1ex>[d]^{a_{3}}
&
&
&
\text{$a_{i}b_{i} = b_{i}a_{i} = 0$, $i = 1, \dots, 4$,}
\\
&
\bullet^{0}
\ar@<1ex>[u]^{a_{1}}
\ar@<1ex>[r]^{b_{4}}
&
\bullet^{3}
\ar@<1ex>[u]^{b_{3}}
\ar@<1ex>[l]^{a_{4}}
&
&
&
\text{$a_{i}a_{j} = b_{l}b_{m} = 0$, when the composition makes sense,}
\\
\G
&
\bullet^{1}
\ar@<1ex>[r]^{d}
&
\bullet^{0}
\ar@<1ex>[l]^{c}
\ar@<1ex>[r]^{a}
&
\bullet^{2}
\ar@<1ex>[l]^{b}
&
&
\text{$ba = dc$; $ab$ and $cd$ nilpotent,}
\\
\H^{m} (m \in \NN)
&
\bullet^{0}
\ar@<1ex>[r]^{a_{1}}
&
\bullet^{1}
\ar@<1ex>[l]^{b_{1}}
\ar@{..}[r]
&
\bullet^{m-1}
\ar@<1ex>[r]^{a_{m}}
&
\bullet^{m}
\ar@<1ex>[l]^{b_{m}}
&
\text{$a_{i}b_{i} = b_{i}a_{i} = 0$, $i = 1, \dots, m$,}
\\
\I^{m} (m \in \NN)
&
\bullet^{1}
\ar@<1ex>[r]^{a_{1}}
\ar@<1ex>[d]^{b_{m}}
&
\bullet^{2}
\ar@<1ex>[l]^{b_{1}}
\ar@<1ex>[d]^{a_{2}}
&
&
&
\text{$a_{i}b_{i} = b_{i}a_{i} = 0$, $i = 1, \dots, m$.}
\\
&
\bullet^{m}
\ar@<1ex>[u]^{a_{m}}
\ar@{..}[r]
&
\bullet^{3}
\ar@<1ex>[u]^{b_{2}}
&
&
&
}
\]
\caption[margin=150pt]{List of tame categories considered in \cite{BB1}. 
If a set of paths is said to be nilpotent, this means that we are considering the inverse limit of the categories such that 
this set of paths is nilpotent of finite nilpotency order.}
\label{fig:tame}
\end{figure}
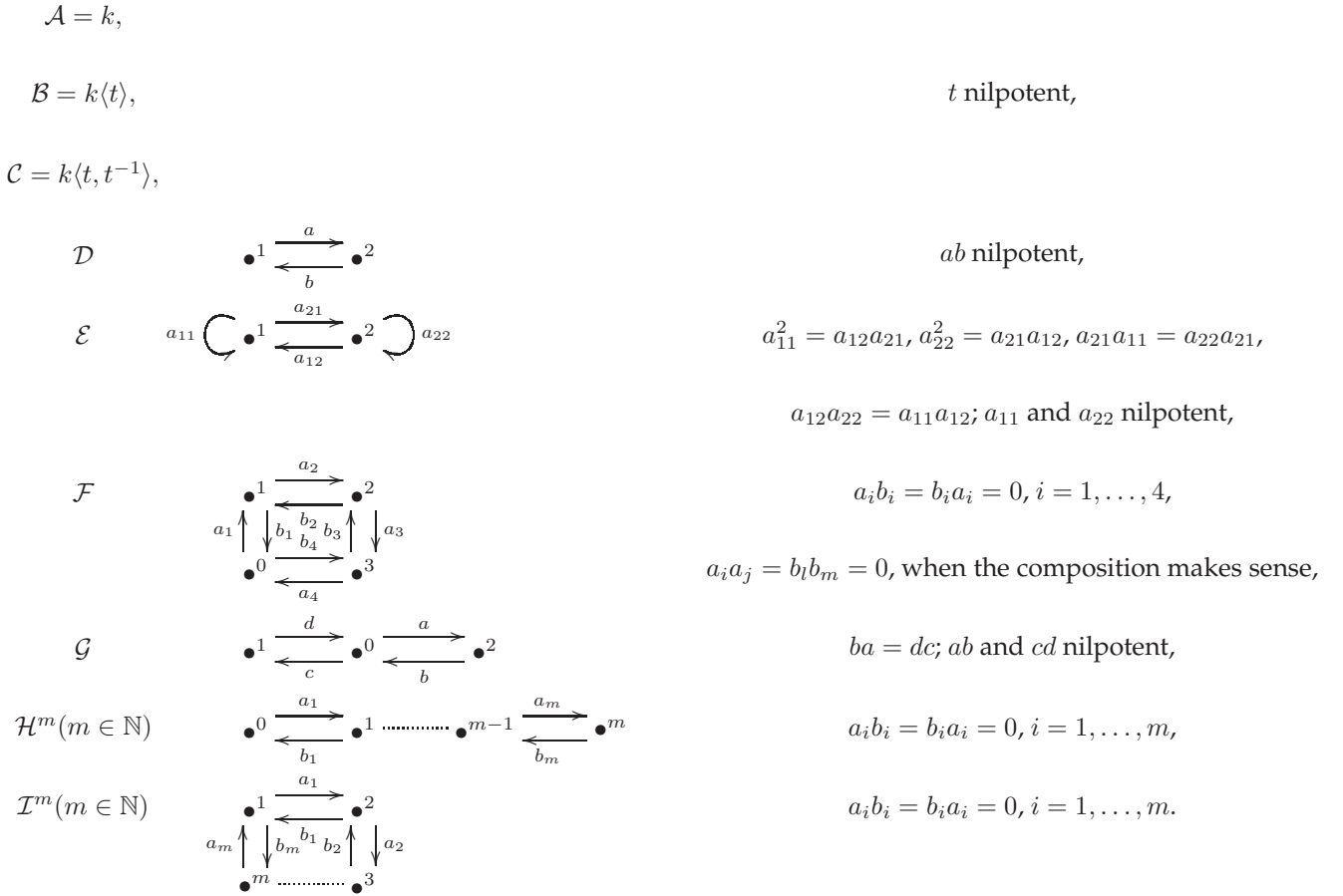
\end{center}
On the other hand, they define the following collection of modules over some of the previous categories,
which verify the properties stated in the next proposition. 
For a given $k$-linear category $\C$, we denote by $\Ind(\C)$ the set of isomorphism classes of indecomposable finite dimensional $\C$-modules. 
\begin{prop}
\begin{itemize}
\item[(i)]
The indecomposable $\B$-modules are given by
\[     \Ind (\B) = \{ \B_{n} : n \in \NN \},     \]
where $\B_{n} = k[t]/(t^{n})$, for $n \in \NN$.

\item[(ii)]
The indecomposable $\C$-modules are given by
\[     \Ind (\C) = \{ \C_{f,n} : n \in \NN, f \in \mathrm{Irr}(k[t]) \setminus \{ t \} \}.     \]
where $\C_{f,n} = k[t]/(f^{n})$, for $n \in \NN$ and $f \in \mathrm{Irr}(k[t]) \setminus \{ t \}$.

\item[(iii)]
The indecomposable $\D$-modules are given by
\[     \Ind (\D) = \{ \D_{n,i} : n \in \NN, i \in \{ 1 , 2 \} \},     \]
where $\D_{n,i}$ ($n \in \NN$, $i = 1,2$) are defined in the $k$-vector space generated by $e_{1}, \dots, e_{n}$, such that
$\D_{n,1}(1)$, $\D_{n,2}(2)$ and $\D_{n,1}(2)$, $\D_{n,2}(1)$ contain the vector $e_{j}$ with odd and even indices, respectively.
The action is given by $u e_{j} = e_{j+1}$ and $v e_{j} = e_{j+1}$ ($e_{n+1}=0$).

\item[(iv)]
The description of the indecomposable $\I^{m}$-modules is more involved and we refer to \cite{BB1}, Sec. 3.2.
We denote by $W$ the free monoid generated by two letters $a$ and $b$ and, given $m \in \NN$, by $\Omega_{m}$
the set of equivalence classes of all non-periodic words.

Given $m \in \NN$,
\[     \Ind (\I^{m}) = \{ \I^{m}_{j,w}, \I^{m}_{z,f} : z \in \Omega_{m}, w \in W, f \in \Ind(k[t]), j = 0, \dots, m - 1 \}.     \]

\item[(v)]
Since the category $\F$ is a quotient of $\I^{4}$, we consider the following $\F$-modules
$\F_{1,x-1} = \I^{4}_{aabb,x-1}$, $\F_{2,x-1} = \I^{4}_{bbaa,x-1}$, $\F_{1,f} = \I^{4}_{abab,f}$, $\F_{2,f} = \I^{4}_{baba,f}$ and
$\F_{j,w} = \I^{4}_{j,w}$, for $f \in \Ind (k[t]), j = 0 , \dots , m - 1$, $w=b^{p}(ab)^{q}a^{r}$, $q \in \NN_{0}$, $p,r \in \{ 0 , 1 \}$.
In this case, the class of indecomposable $\F$-modules is
\[     \Ind (\F) = \{ \F_{1,x-1}, \F_{2,x-1}, \F_{1,f}, \F_{2,f}, \F_{j,w} : f \in \Ind (k[t]), j = 0 , \dots , m - 1, w = b^{p}(ab)^{q}a^{r},
       q \in \NN_{0}, p, r \in \{ 0 , 1 \} \}.
\]

\item[(vi)]
Using the functor $\H^{m} \rightarrow \I^{m+2}$, given by $i \mapsto i$, $a_{i} \mapsto a_{i}$ y $b_{i} \mapsto b_{i}$, we are able to define
the family of $\H^{m}$-modules $\H^{m}_{j,w} = \I_{j,w}^{m+2}$, for $w \in W$, $|w| \leq m - j$, $j = 0, \dots, m - 1$.
The indecomposable $\H^{m}$-modules are given by
\[     \Ind (\H^{m}) = \{ \H^{m}_{j,w} : w \in W, |w| \leq m - j, j = 0, \dots, m - 1 \}.     \]
\end{itemize}
\end{prop}
\noindent\textbf{Proof.}
Cf. \cite{BB1}, Sec. 3.1. to 3.4.
\qed

We may now state the results describing the subcategories of weight modules and generalized weight modules over the Weyl algebras $A_{r}(k)$.
\begin{teo}
Let us consider the first Weyl algebra $A_{1}(k)$.
In this case, the subcategory $\W(\O)$ is tame for each orbit $\O$.
Moreover,
\begin{itemize}
\item[(i)] if $\O \neq \ZZ$, then $\W(\O) \simeq {}_{\A}\Mod$ and $\Ind(\W(\O)) = \{ \A_{\O} \}$.

\item[(ii)] if $\O = \ZZ$, then $\W(\O) \simeq {}_{\H^{1}}\Mod$ and
$\Ind(\W(\O)) = \{ \H^{1}_{\O,j,w} \}$.
\end{itemize}

On the other hand, the category $\GW(\O)$ is tame since $\mathrm{char}(k) = 0$.
Besides, in case $\O \neq \ZZ$, $\GW(\O) \simeq {}_{\B}\Mod$ and $\Ind(\GW(\O)) = \{ \B_{\O, n} : n \in \NN \}$.
Finally, if $\O = \ZZ$, then $\GW(\O) \simeq {}_{\D}\Mod$ and $\Ind(\GW(\O)) = \{ \D_{\O, n, i} : n \in \NN, i= 1,2 \}$.
\end{teo}
\noindent\textbf{Proof.}
Cf. \cite{BB1}, Sec. 4.1.
\qed

\begin{teo}
Let us consider the second Weyl algebra $A_{2}(k)$.
Then, $\W(\O)$ is tame for each orbit $\O = \O_{1} \times \O_{2}$.
Moreover,
\begin{itemize}
\item[(i)] if $\O_{1} \neq \ZZ$ and $\O_{2} \neq \ZZ$, $\W(\O) \simeq {}_{\A}\Mod$ and $\Ind(\W(\O)) = \{ \A_{\O} \}$.
\item[(ii)] if $\O_{i} = \ZZ$ and $\O_{j} \neq \ZZ$ ($i,j \in \{1,2\}$), then $\W(\O) \simeq {}_{\H^{1}}\Mod$ and $\Ind(\W(\O))
= \{ \H^{1}_{\O,j,w} \}$.
\item[(iii)] if $\O_{1} =  \O_{2} = \ZZ$, $\W(\O) \simeq {}_{\F}\Mod$ and
$\Ind(\W(\O)) = \{ \F_{\O,1,x-1}, \F_{\O,2,x-1}, \F_{\O,1,f}, \F_{\O,2,f}, \F_{\O,j,w} \}$.
\end{itemize}

The subcategory $\GW(\O)$ is wild for each orbit $\O$.
\end{teo}
\noindent\textbf{Proof.}
Cf. \cite{BB1}, Sec. 4.2.
\qed

\begin{teo}
Let us now consider the Weyl algebra $A_{n}(k)$, $n \geq 3$.
We shall write each orbit $\O = \O_{1} \times \dots \times \O_{n}$.
\begin{itemize}
\item[(i)] If $\O_{i} \neq \ZZ$ ($i = 1, \dots, n$), then $\W(\O) \simeq {}_{\A}\Mod$ and $\Ind(\W(\O)) = \{ \A_{\O} \}$.

\item[(ii)] If the number of non-degenerate orbits $\O_{i}$ is $n-1$, then $\W(\O) \simeq {}_{\H^{1}}\Mod$ and $\Ind(\W(\O)) = \{ \H^{1}_{\O,j,w} \}$.

\item[(iii)] If the number of no-degenerate orbits $\O_{i}$ is $n-2$, then $\W(\O) \simeq {}_{\F}\Mod$ and
\[     \Ind(\W(\O)) = \{ \F_{\O,1,x-1}, \F_{\O,2,x-1}, \F_{\O,1,f}, \F_{\O,2,f}, \F_{\O,j,w} \}.     \]

\item[(iv)] If the number of non-degenerate orbits $\O_{i}$ is less than $n-2$, then $\W(\O)$ es wild.
\end{itemize}
\end{teo}
\noindent\textbf{Proof.}
Cf. \cite{BB1}, Sec. 4.3.
\qed

Finally, we shall give the proof of Theorem \ref{teo:yangmillsweyl}.

\noindent\textbf{Proof of Theorem \ref{teo:yangmillsweyl}.}
We have studied in detail cases $r = 1, 2, 3, 4$.
We shall then restrict ourselves to the case $r \geq 5$.

We know that $\ym(3) = V(3) \oplus \tym(3)$ as $k$-modules.
Also, considered as a graded Lie algebra with the special grading, the Lie ideal $\tym(3)$ is a free graded Lie algebra
(concentrated in even degrees) generated by a graded vector space (concentrated in even degrees) $W(3)$, that is,
\[     \tym(3) \simeq \f_{gr}(W(3))  \hskip 5mm  \text{and} \hskip 5mm  W(3) = \bigoplus_{l \in \NN} W(3)_{2 l + 2},     \]
where $\dim_{k}(W(3)_{2 l + 2}) = 2 l + 1$ (cf. Eq. \eqref{eq:W(3)}).
In fact, we shall choose a basis for the first two homogeneous subspaces of $W(3)$ as we did previously.
For $W(3)_{4}$ we fix the basis $\{x_{12}, x_{13}, x_{23} \}$ and for $W(3)_{6}$ we fix the basis
$\{ x_{112}, x_{113}, x_{221}, x_{123}, x_{312} \}$.

By Proposition \ref{prop:conmutacionfuntores},
\[     \U_{gr}(\tym(3)) \simeq \U_{gr}(\f_{gr}(W(3))) \simeq T_{gr}(W(3)),     \]
and hence,
\[     T(\mathcal{O}(W(3))) = \mathcal{O}(T_{gr}(W(3))) \simeq \mathcal{O}(\U_{gr}(\tym(3))) = \U(\mathcal{O}(\tym(3))),     \]
where $\mathcal{O}$ denotes the corresponding forgetful functors. 

Having a morphism of $k$-algebras from $\U(\mathcal{O}(\tym(3)))$ to $A_{m}(k)$ is the same as having a morphism of $k$-algebras from
$T(\mathcal{O}(W(3)))$ to $A_{m}(k)$, which in turn is the same as having a morphism of $k$-vector spaces from $\mathcal{O}(W(3))$ to $A_{m}(k)$.
The morphism of algebras will be surjective if the image of the corresponding morphism of vector spaces contains the generators as an algebra
of $A_{m}(k)$, denoted by $p_{1}, \dots, p_{m}, q_{1}, \dots, q_{m}$ as usual.

From now on, we shall exclusively work in the non-graded case (for algebras), so we will omit the forgetful functor $\mathcal{O}$ without confusion.
However, we will keep the canonical grading of the Yang-Mills algebra $\ym(3)$.

Let us suppose $m \geq 2$.
We consider the morphism of $k$-vector spaces
\[     \phi : W(3) \rightarrow A_{m}(k)     \]
such that $\phi(W(3)_{4}) = \{0\}$, $\phi(x_{112}) = p_{1}$, $\phi(x_{221}) = p_{2}$, $\phi(x_{123}) = q_{1}$, $\phi(x_{113}) = q_{2}$,
$\phi(x_{312}) = 0$,
and such that for each generator of the Weyl algebra $p_{i}$ and $q_{i}$ ($i \geq 3$), there exist homogeneous elements of degree greater than
$6$, $w_{i}, w'_{i} \in W(3)$ such that $\phi(w_{i}) = p_{i}$ and $\phi(w'_{i}) = q_{i}$.
This last condition is easily verified, taking into account that $W(3)$ is infinite dimensional.
Of course this means that there are a lot of choices for this morphism.

Let $d_{i}$ and $d'_{i}$ be the degrees of $w_{i}$ and $w'_{i}$, respectively.
Let $j$ be the maximum of the degrees $d_{i}$ and $d'_{i}$, and let $l = 2 j + 1$.
The morphism $\phi$ induces a unique surjective homomorphism $\Phi : \U(\tym(3)) \twoheadrightarrow A_{m}(k)$, equivalent to the
homomorphism of Lie algebras
\[     \tym(3) \rightarrow \mathrm{Lie}(A_{m}(k)),     \]
where $\mathrm{Lie}(\place) : {}_{k}\Alg \rightarrow {}_{k}\LieAlg$ denotes the functor that associates to every associative algebra 
the Lie algebra with the same underlying vector space and Lie bracket given by the commutator of the algebra. 
The latter morphism may be factorized in the following way
\[     \tym(3) \rightarrow \tym(3)/\mathcal{C}^{l}(\ym(3)) \rightarrow \mathrm{Lie}(A_{m}(k)),     \]
where the first morphism is the canonical projection.
Hence, the map $\Phi$ may be factorized as 
\[     \U(\tym(3)) \twoheadrightarrow \U(\tym(3)/\mathcal{C}^{l}(\ym(3))) \twoheadrightarrow A_{m}(k).     \]
We have thus obtained a surjective homomorphism of $k$-algebras
\[     \Psi : \U(\tym(3)/\mathcal{C}^{l}(\ym(3))) \twoheadrightarrow A_{m}(k),     \]
where the Lie algebra $\tym(3)/\mathcal{C}^{l}(\ym(3))$ is nilpotent.
Moreover, it is a nilpotent ideal of the (nilpotent) Lie algebra $\ym(3)/\mathcal{C}^{l}(\ym(3))$.
We have, as $k$-modules,
\[     \ym(3)/\mathcal{C}^{l}(\ym(3)) = V(3) \oplus \tym(3)/\mathcal{C}^{l}(\ym(3)).     \]

Let $I$ be the kernel of $\Psi$ in $\U(\tym(3)/\mathcal{C}^{l}(\ym(3)))$.
Taking into account that the quotient of the universal enveloping algebra $\U(\tym(3)/\mathcal{C}^{l}(\ym(3)))$ by $I$ is a Weyl algebra which is simple,
$I$ is a maximal two-sided ideal, and then, there exists a linear functional
$f \in (\tym(3)/\mathcal{C}^{l}(\ym(3)))^{*}$ such that $I = I(f)$.
We fix a standard polarization $\h_{f}$ for $f$, i.e., a polarization constructed from a flag of ideals of $\tym(3)/\mathcal{C}^{l}(\ym(3))$.
Let $\bar{f} \in (\ym(3)/\mathcal{C}^{l}(\ym(3)))^{*}$ be any extension of $f$, and let $\h_{\bar{f}}$ be a standard polarization for $\bar{f}$
given by extending the flag of $\tym(3)/\mathcal{C}^{l}(\ym(3))$ to $\ym(3)/\mathcal{C}^{l}(\ym(3))$, i.e., if the flag is:
\[     \tym(3)/\mathcal{C}^{l}(\ym(3)) \subset \g_{1} \subset \g_{2} \subset \g_{3} = \ym(3)/\mathcal{C}^{l}(\ym(3)),     \]
and we denote $\bar{f}_{i}$ the restriction of $\bar{f}$ to $\g_{i}$ ($i = 1,2,3$), then
\[     \h_{\bar{f}} = \h_{f} + \g_{1}^{\bar{f}_{1}} + \g_{2}^{\bar{f}_{2}} + \g_{3}^{\bar{f}_{3}}.     \]

Let $\g$ be a finite dimensional Lie algebra and $\h$ a subalgebra.
Given $I \triangleleft \U(\h) \subset \U(\g)$ a two sided ideal in an enveloping algebra of $\h$, consider
\[     \mathfrak{st}(I, \g) = \{ x \in \g : [x , I] \subset I \}.     \]

Since $\tym(3)/\mathcal{C}^{l}(\ym(3))$ is an ideal of the nilpotent Lie algebra $\ym(3)/\mathcal{C}^{l}(\ym(3))$, according to Proposition 6.2.8 of
\cite{Dix1}
\[     \mathfrak{st}(I(f), \ym(3)/\mathcal{C}^{l}(\ym(3))) = \tym(3)/\mathcal{C}^{l}(\ym(3)) + \g',     \]
where
\[     \g' = \{ x \in \ym(3)/\mathcal{C}^{l}(\ym(3)) : f([x , \tym(3)/\mathcal{C}^{l}(\ym(3))]) = 0 \}.     \]

For our ideal, we get immediately that $\bar{x}_{12}, \bar{x}_{13}, \bar{x}_{23} \in I$, but $\bar{x}_{112}, \bar{x}_{221}, \bar{x}_{123}$ do not
belong to $I$, since $\Psi(\bar{x}_{112}) = p_{1}$, $\Psi(\bar{x}_{221}) = p_{2}$ and $\Psi(\bar{x}_{123}) = q_{1}$.

Let $x \in \ym(3)/\mathcal{C}^{l}(\ym(3))$, then $x = x' + y$, where
\[     x' = \sum_{i=1}^{3} c_{i} \bar{x}_{i} \in V(3),     \]
and $y \in \tym(3)/\mathcal{C}^{l}(\ym(3))$.
Since $[y,I(f)] \subset I(f)$, this implies that $x \in \mathfrak{st}(I(f), \ym(3)/\mathcal{C}^{l}(\ym(3)))$ if and only if
$[x' , I(f)] \subset I(f)$.
Explicitly,
\[     [x' , \bar{x}_{12}] = \sum_{i=1}^{3} c_{i} [\bar{x}_{i} , \bar{x} _{12}] = c_{1} \bar{x}_{112} - c_{2} \bar{x}_{221} - c_{3} \bar{x}_{123}.     \]
If $[x' , \bar{x}_{12}] \in I$, then $\Psi([x' , \bar{x}_{12}]) = 0$, or,
\[     c_{1} p_{1} - c_{2} p_{2} - c_{3} q_{1} = 0,     \]
but the generators of $A_{m}(k)$ are linearly independent, so $c_{1} = c_{2} = c_{3} = 0$, which gives $x' = 0$, implying
$\mathfrak{st}(I(f), \ym(3)/\mathcal{C}^{l}(\ym(3))) = \tym(3)/\mathcal{C}^{l}(\ym(3))$.
Hence $g' \subset \tym(3)/\mathcal{C}^{l}(\ym(3))$.
As a consequence, we obtain the inclusion $\g_{i}^{\bar{f}_{i}} \subset \tym(3)/\mathcal{C}^{l}(\ym(3))$, whence
$\h_{\bar{f}} \subset \tym(3)/\mathcal{C}^{l}(\ym(3))$.
By maximality of $\h_{f}$ in $\tym(3)/\mathcal{C}^{l}(\ym(3))$, we find that $\h_{\bar{f}} \subset \h_{f}$.
The other inclusion is even simpler, so $\h_{\bar{f}} = \h_{f}$.

Finally, the weight of the ideal $I(\bar{f})$ is equal to
\begin{align*}
     \dim_{k}((\ym(3)/\mathcal{C}^{l}(\ym(3)))/\h_{\bar{f}}) &= \dim_{k}((\ym(3)/\mathcal{C}^{l}(\ym(3)))/\h_{f})
     \\
     &= \dim_{k}((\tym(3)/\mathcal{C}^{l}(\ym(3)))/\h_{\bar{f}}) + 3 = m + 3.
\end{align*}

We have then proved that $A_{r}(k)$ is a quotient of $\U(\ym(3))$, for any $r \geq 5$, and, as a consequence, for any $r \in \NN$.
\qed

\begin{obs}
The previous proof does not work for $\ym(2)$, since in this case $\tym(2) \simeq \f_{gr}(W(2))$ with $\dim(W(2)) = 1$.
\end{obs}

\begin{defi}
Let $A$ be an associative or Lie $k$-algebra, and let $\R \subset {}_{A}\Mod$ be a full subcategory of the category of (left) modules over over $A$.
We shall say that $\R$ \emph{separates points} of $A$ if for all $a \in A$, $a \neq 0$, there exists $M \in \R$ such that $a$ acts on $M$ as a non
null morphism.
\end{defi}

By Proposition 3.1.15 of \cite{Dix1}, we know that an ideal $I$ in the enveloping algebra of a finite dimensional Lie algebra $\g$ is
semiprime if and only if it is the intersection of primitive ideals.
Since $\U(\g)$ is a domain, the null ideal $\{ 0 \}$ is completely prime, and hence semiprime (cf. \cite{Dix1}, 3.1.6).
As a consequence, $\{ 0 \}$ is an intersection of two-sided primitive ideals.
A fortiori, the intersection of all primitive two-sided ideals in $\U(\g)$ is null, that is, the Jacobson radical of the
enveloping algebra of $\g$ is null, i.e. $J(\U(\g)) = \{ 0 \}$.
Analogously, since every maximal two-sided ideal is primitive (cf. \cite{Dix1}, 3.1.6), the intersection of all maximal two-sided ideals is null.

Let $\W(n)$ be the full subcategory of ${}_{\U(\ym(n))}\Mod$ consisting of all modules $M$ satisfying the following property:
there exist $r, l \in \NN$ such that the action of the Yang-Mills algebra may be factorized as follows
\[     \phi : \U(\ym(n)) \twoheadrightarrow \U(\ym(n)/\C^{l}(\ym(n))) \rightarrow A_{r}(k) \rightarrow \End_{k}(M).     \]

\begin{prop}
If $n \geq 3$, the category $\W(n)$ separates points of $\U(\ym(n))$.
\end{prop}
\noindent\textbf{Proof.}
Let $x \in \U(\ym(n))$ be a non zero element.
There exists then $l \in \NN$ such that $\pi_{l}(x) \in \U(\ym(n)/\C^{l}(\ym(n)))$ is non zero (cf. Lemma \ref{lema:yangmillsnulo}).
Since $\ym(n)/\C^{l}(\ym(n))$ is a nilpotent finite dimensional Lie algebra, the intersection of all maximal two-sided ideals is zero,
so there exists a maximal two-sided ideal $J \triangleleft \U(\ym(n)/\C^{l}(\ym(n)))$ such that $\pi_{l}(x) \notin J$.
On the other hand, $\U(\ym(n)/\C^{l}(\ym(n)))/J \simeq A_{r}(k)$, for some $r \in \NN$ (cf. \cite{Dix1}, 4.5.8 and Thm. 4.7.9).
Taking into account that the inverse image of a maximal two-sided ideal by a surjective $k$-algebra homomorphism is maximal, $I = \pi_{l}^{-1}(J)$
is a maximal two-sided ideal in $\U(\ym(n))$ and $\U(\ym(n))/I \simeq \U(\ym(n)/\C^{l}(\ym(n)))/J \simeq A_{r}(k)$.
If $x \in I$, then $\pi_{l}(x) \in \pi_{l} (I) \subset J$, thus $x \notin I$.

Let $M$ be a (left) module over $A_{r}(k)$, such that the image of $x$ under the previous isomorphisms is not zero (take for instance $A_{r}(k)$).
Hence, the previous isomorphisms induce a structure of $\U(\ym(n))$-module over $M$, such that $x$ does not act as the zero endomorphism on $M$.
The proposition is then proved.
\qed

\begin{rem}
\label{obs:separa}
Although the subcategory $\W(n)$ separates points of the Yang-Mills algebra $\YM(n)$, it does not satisfy that every element 
in the category ${}_{\YM(n)}\Mod$ is isomorphic to an element in $\W(n)$. 
This is a consequence of the following fact: every object of $\W(n)$ has finite Gelfand-Kirillov dimension, for, if $M \in \W(n)$ is an induced module of a module
over $A_{r}(k)$, then (cf. \cite{McR1}, Prop. 1.15, (ii); Prop. 3.2, (iii), (v))
\[     \text{$\mathrm{GK}$-$\dim_{\YM(n)}(M)$} = \text{$\mathrm{GK}$-$\dim_{A_{r}(k)}(M)$} \leq  \text{$\mathrm{GK}$-$\dim_{A_{r}(k)}(A_{r}(k))$} = 2r.     \]
On the other hand, $\YM(n)$ having infinite Gelfand-Kirillov dimension since it has exponential growth (cf. \cite{CD1}), there exist $\YM(n)$-modules
(e.g. the regular module $\YM(n)$ in itself) which do not belong to $\W(n)$.
\end{rem}



\end{document}